\documentclass[final,leqno]{siamltex}

\usepackage{amsfonts,amssymb} 
\usepackage{amsmath} 
\usepackage{color} 
\usepackage{graphicx} 
\usepackage{color}
\usepackage{epsfig,mathrsfs} 
\usepackage[bf,SL,BF]{subfigure} 
\usepackage{fancyhdr} 
\usepackage{CJK} 
\usepackage{caption} 
\usepackage{wrapfig} 
\usepackage{cases} 
\usepackage{bm}
\usepackage{algorithm}
\usepackage{algorithmic}
\usepackage{enumerate}
\usepackage{multirow}
\usepackage{threeparttable}
\usepackage{txfonts}
\usepackage[colorlinks]{hyperref}
\usepackage{hypernat}
\newtheorem{remark}{Remark}[section] 

\makeatletter
\newcommand{\rmnum}[1]{\romannumeral #1}
\newcommand{\Rmnum}[1]{\expandafter\@slowromancap\romannumeral #1@}
\makeatother

\title{A new extrapolation cascadic multigrid method for 3D elliptic boundary value problems on rectangular domains\thanks{This research was supported by the Natural Science Foundation of China (Nos. 41474103, 41204082, 11402174 and 11301176), the National
High Technology Research and Development Program of China (No. 2014AA06A602) and
the Natural Science Foundation of Hunan Province of China (No. 2015JJ3148).}}

\author{Kejia Pan\footnotemark[2]\ \footnotemark[3]
\and Dongdong He\footnotemark[4]
\and Hongling Hu\footnotemark[5]}

\begin{document}
\maketitle

\renewcommand{\thefootnote}{\fnsymbol{footnote}}

\footnotetext[2]{School of Mathematics and Statistics, Central South University, Changsha 410083, P.R. China (E-mail: pankejia@hotmail.com)}
\footnotetext[3]{Department of Mathematics and Statistics, University of North Carolina at Charlotte, Charlotte, NC 28223, USA}
\footnotetext[4]{School of Aerospace Engineering and Applied Mechanics, Tongji University, Shanghai 200092, P.R. China (E-mail: dongdonghe@tongji.edu.cn)}
\footnotetext[5]{Corresponding author (hhling625@163.com). College of Mathematics and Computer Science, Key Laboratory of High Performance Computing and Stochastic Information Processing  (Ministry of Education of China), Hunan Normal University, Changsha 410081, P.R. China}

\renewcommand{\thefootnote}{\arabic{footnote}}

\begin{abstract}
 In this paper, we develop a new extrapolation cascadic multigrid method ($\textrm{ECMG}_{jcg}$), which makes it possible to solve 3D elliptic boundary value problems on rectangular domains of over 100 million unknowns on a desktop computer in minutes.
  First, by combining Richardson extrapolation and tri-quadratic Serendipity interpolation techniques, we introduce a new extrapolation formula to provide a
 good initial guess for the iterative solution on the next finer grid, which is a third order approximation to the finite element (FE) solution.
  And the resulting large sparse linear system from the FE discretization is then solved by the Jacobi-preconditioned Conjugate Gradient (JCG) method.
  Additionally, instead of performing a fixed number of iterations as used in the most of cascadic multigrid method (CMG) literature, a relative residual stopping criterion is used in our iterative solvers, which enables us to obtain conveniently the numerical solution with the desired accuracy.
  Moreover, a simple Richardson extrapolation is used to cheaply get a fourth order accurate
solution on the entire fine grid from two second order accurate solutions on two different scale grids.
 Test results from three different problems  with smooth and singular solutions are reported to show that $\textrm{ECMG}_{jcg}$ has much better efficiency compared to  the classical V-cycle and W-cycle multigrid methods.
  Since the initial guess for the iterative solution is a quite good approximation to the FE solution, numerical results show that only few number of iterations are required on the finest grid for  $\textrm{ECMG}_{jcg}$ with an appropriate tolerance of the relative residual to achieve full second order accuracy, which is particularly important when solving large systems of equations and can greatly reduce the computational cost. It should be pointed out that when the tolerance becomes more cruel,  $\textrm{ECMG}_{jcg}$ still needs only few iterations to obtain fourth order extrapolated solution on each grid, except on the finest grid. Finally, we present the reason why our ECMG algorithms are so highly efficient for solving these elliptic problems.
\end{abstract}



\pagestyle{myheadings}
\thispagestyle{plain}
\markboth{K. J. PAN AND D. D. HE AND H. L. HU}{A new ECMG method for 3D elliptic boundary value problems}

\section{Introduction}
Elliptic boundary value problems arise in many areas of geophysical fluid dynamics. Consider the following model problem:
 \begin{equation}\label{bvp}
\left\{ \begin{aligned}
         -\nabla\cdot\big(\beta(\mathbf{x})\nabla u\big)&=&f(\mathbf{x})\quad  &\textrm{in } \Omega,\\
          u&=&g_D(\mathbf{x})\quad &\textrm{on } \Gamma_D,\\
          \alpha(\mathbf{x}) u+\beta(\mathbf{x}) \frac{\partial u}{\partial n}&=&g_R(\mathbf{x})\quad &\textrm{on } \Gamma_R,
        \end{aligned} \right.
\end{equation}
where $\alpha$ and $\beta$ are piecewise smooth functions on $\overline{\Omega}$ and $0<\beta_{min}\leq \beta\leq \beta_{max}$ for every
$\mathbf{x}\in \Omega$, $\mathbf{n}$ is the outward unit normal to $\partial\Omega$, $f:\Omega\rightarrow\Re, g_D:\Gamma_D\rightarrow\Re$ and $g_R:\Gamma_R\rightarrow\Re$  are assigned functions. Here
$\Omega$ is some bounded rectangular domain in $R^3$ with Dirichlet boundary $\Gamma_D$ and
Robin boundary $\Gamma_R$. It is well known that the Neumann boundary condition corresponds to the extreme
case, namely $\alpha = 0$.

The elliptic boundary value problem can be approximated using different numerical techniques, such as finite difference (FD) and finite element (FE) methods.
The resulting linear system can be solved efficiently using direct solver for problems with less than one million unknowns.
However, for 3D problems, even only a few hundreds grid points in each coordinate direction already leads  to
a system with millions of unknowns, one has to resort to an iterative method.
For example, when solving direct current (DC) resistivity and electromagnetic modelling problems arising in geophysical applications,
the grids designed to approximate huge realistic 3D geologies are
usually enormous in order to represent correctly complex structures. Consequently,
it is normally necessary to solve hundreds of millions unknowns in the forward problem.
In addition, the forward problem has to be solved many times in the inversion of geophysical data \cite{Avdeev2005,Newman2014,Koldan2014}.
Therefore, it is critical that the 3D problem is solved very efficiently.
The Multigrid technique~\cite{Briggs2000, Trottenberg2001} is one of the most efficient strategies to solve the large linear system
from discretized elliptic differential equations.
 The classical MG methods have been successfully applied to solve the Poisson equations~\cite{Schaffer1984, Gupta1997a,Othman1999, Zhang1998, Zhang2002, Wang2009, Ge2010}, Helmholtz
equation~\cite{Erlangga2006, Elman2001, Erlangga20062, Kim2002, Riyanti2007} and  convection-diffusion equations~\cite{Zeeuw1995,Gupta1997,Gupta2000,Zhang20022,Ge2011}. However, traditional
MG methods (both geometric and algebraic) have to cycle between
coarse and fine grids in order to accelerate the rate of convergence.
Therefore, MG methods are difficult to implement in programming
language.

The Cascadic multigrid (CMG) method proposed by Deuflhard and Bornemann in \cite{Bornemann1996} is a simpler multilevel method
without coarse-grid correction. The CMG method uses CG solvers as the basic iteration methods on successively
refined grids where the initial guesses are the linear interpolations
of the approximate solutions on the previous grids. Nevertheless,
the CMG method has the same optimal property compared to MG
methods. Namely, the algorithm converges with a rate that is independent
of the grid sizes and the numbers of grids levels \cite{Bornemann1996,Shaidurov1996}. In 2008, an extrapolation cascadic
multigrid (ECMG) method was presented by Chen et al. in \cite{Chen2008} for
solving the second-order elliptic boundary value problems.
This method proposes a new extrapolation formula to provide a better initial guess for the iterative solution on the next finer grid,
which improves the convergence rate of the original CMG algorithm.
However, as far as we know, the ECMG algorithm has mainly been used for solving the 2D elliptic
boundary value problems in existing literature. But it is of more importance to solve the 3D problems efficiently and accurately arising in many engineering areas, such as geophysical exploration \cite{Avdeev2005,Newman2014}. And it is nontrivial to extend the ECMG method from 2D to 3D.

 In this paper, we develop a new extrapolation cascadic multigrid method ($\textrm{ECMG}_{jcg}$) for solving the 3D elliptic boundary value problems on rectangular domains.  In our approach, the computational domain is discretized by a regular hexahedral grid,
 and a linear FE method is used to discretize the 3D elliptic problem.
By extrapolating  the numerical solutions on a coarse grid and a fine grid (with half the mesh size) to obtain a good initial guess of the iterative solution
 on the 8 vertices and 12 edge-midpoints of each coarse hexahedral element, which consists of 64 connected small hexahedral elements of the next finer grid (two times refined grid with one-fourth the mesh size), and then by using tri-quadratic Serendipity interpolation for the above 20 nodes to get the initial guesses on the other 105 $(5^3-20)$ nodes of such coarse hexahedral element,
  we are able to obtain a third order approximation to the FE solution on the next finer grid.
  And the resulting large sparse linear system from the FE discretization is then solved by the Jacobi-preconditioned Conjugate Gradient (JCG) solver using the obtained initial guess.
    Additionally, a tolerance related to relative residual is introduced in the JCG solver.
     Moreover, a simple Richardson extrapolation is used to obtain cheaply a fourth order accurate solution on the entire fine grid from two second order accurate solutions on two different scale grids (current fine and previous coarser grids).
     Finally, our method has been used to solve 3D elliptic problems with more than 135 million unknowns in about 1 minute on a small server with 32GB RAM installed. Compared to the existing ECMG method, our new ECMG approach is advantageous due to the following reasons:
   \begin{itemize}
  \item Instead of performing a fixed number of iterations as used in the usual CMG methods, a relative residual tolerance is used in our iterative solvers, which enables us to not only avoid the difficulty of choosing the number of iterations on each grid,
  but also allows us to obtain conveniently the numerical solution with the desired accuracy.

 \item The second important point is to employ JCG method as MG smoother, which requires only a little additional work than CG, while JCG smoother typically yields better convergence properties than CG smoother.

  \item A fourth order extrapolated solution on the entire fine grid is constructed to greatly enhance the accuracy of the numerical solution.

  \item we can clearly explain why only few number of iterations are required on the finest grid to achieve full second order accuracy for our method through defining a ratio.
\end{itemize}

The rest of the paper is organized as follows:
Section 2 gives the description of the FE discretization for the 3D elliptic boundary value problem. Section 3 reviews the classical V-cycle and W-cycle MG methods. In Section 4, we first present some extrapolation and interpolation formulas, and then develop a new ECMG method to solve 3D elliptic boundary value problems. Section 5 contains the numerical results to demonstrate the high efficiency and accuracy of the proposed method. And conclusions are given in the final Section.

\section{Finite element discretizations}
In this Section we introduce a linear FE discretization to (\ref{bvp}) which we shall use in the construction of MG methods.
For simplicity, we assume  that eq.(\ref{bvp}) has homogeneous  Dirichlet boundary conditions, namely $g_D\equiv 0$.
We introduce a bilinear form, $a(\cdot, \cdot): H_D^1(\Omega)\times H_D^1(\Omega)\rightarrow R$, defined in the usual way
\begin{equation}
  a(u,v)=\int_{\Omega} \beta \nabla u\cdot \nabla v dx + \int_{\Gamma_R} \alpha uv dx,
\end{equation}
where $H_D^1(\Omega)=\{v\in H^1(\Omega): v|_{\Gamma_D}=0\}$, and $v|_{\Gamma_D}=0$ are in the sense of trace.

Then the weak formulation of the problem (\ref{bvp}) is to find $u\in H_D^1(\Omega)$ such that
\begin{equation}\label{var}
  a(u,v)=f(v), \quad \forall v\in H_D^1(\Omega),
\end{equation}
where
\begin{equation}
  f(v)=\int_{\Omega} fv dx + \int_{\Gamma_R} g_R v ds.
\end{equation}

Assuming that $\Omega$ is partitioned by an uniform hexahedral mesh $T_h$ with characteristic mesh size $h$,
namely $\overline{\Omega}=\bigcup_{\tau\in T_h} \overline{\tau}$. Then a piecewise linear FE space, $V_h\subset H(\Omega)$ can be defined by
\begin{equation}
  V_h=\left\{v\in H_D^1(\Omega):v|_{\tau}\in P_{1}(\tau), \forall \tau\in T_h\right\},
\end{equation}
where $P_1$ denotes the set of linear polynomials.

Let $\{\phi_i\}_{i=1}^N$ be the standard nodal basis functions of FE space $V_h$ and writing the FE solution as the linear combination of the basis functions as:
$u_h=\sum_{i=1}^N u_{h,i}\phi_i$, we can obtain the following algebraic systems for the variational problem (\ref{var}):
\begin{equation}\label{fe_equation}
  A_h u_h = f_h,
\end{equation}
where
\begin{equation}
  A_h = \left(a(\phi_i, \phi_j)\right),\quad u_h = \left(u_{h,i}\right),\quad f_h=\left(f(\phi_i)\right).
\end{equation}
In this paper, we focus on how to solve the FE equation (\ref{fe_equation}) efficiently using a new extrapolation cascadic multigrid method.

\section{Classical multigrid method}
The MG method is very useful in increasing the efficiency of iterative methods (such as Gauss-Seidel or weighted Jacobi) for solving the
large systems of algebraic equations resulting from FE discretizations of elliptic boundary value problems.
 The main idea behind the MG methods is to split the errors into high frequency errors and low frequency errors and
solve them in different subspaces. To be more specific,
the high frequency errors are solved on the fine
grid, while the low frequency components are solved on the coarse grid. A good introductory text on multigrid is the book by Briggs \cite{Briggs2000}; and  a more advanced treatment is given by Trottenberg in \cite{Trottenberg2001}.


\begin{algorithm}[htb]         

\caption{Classical MG Method (Recursive Definition): $u_h$ $\Leftarrow$ MGM($A_h, u_h, f_h$)}             

\label{alg:mg}                  

\begin{algorithmic}[1]                

\IF    {$h==H$(Coarsest grid)}
\STATE $u_H$ $\Leftarrow$ DSOLVE($A_H u_H=f_H$)
\ELSE
  \STATE  Perform $N_{pre}$ pre-smoothing iterations with the current estimate $u_h$
  \STATE $r_h=f_h-A_h u_h$
  \STATE $r_{2h} = R_h^{2h} r_h$
  \STATE $e_{2h}=0$
  \FOR {$i=1$ to $N_{cycles}$}
    \STATE $e_{2h}\Leftarrow$ MGM$(A_{2h}, e_{2h}, r_{2h})$
  \ENDFOR
  \STATE $e_h = P_{2h}^h e_{2h}$
  \STATE $u_h = u_h + e_h$
  \STATE Perform $N_{post}$ post-smoothing iterations with initial guess $u_h$
\ENDIF              

\end{algorithmic}

\end{algorithm}

\begin{algorithm}[htb]         

\caption{ Multigrid solver for the finest grid system $A_h u_h = f_h$}            

\label{alg:fmg}                  

\begin{algorithmic}[1]                

\STATE $u_h =0$
\WHILE {$||A_h u_h-f_h||_2>\epsilon ||f_h||_2$}
    \STATE $u_h \Leftarrow$ MGM$(A_h, u_h, f_h)$
\ENDWHILE

\end{algorithmic}

\end{algorithm}

The MG method can be defined recursively as follows.
On the coarsest grid, a direct solver DSOLVE is used since the size of the linear system is small, see line 2 in the Algorithm \ref{alg:mg}.
The procedure MGM$(A_h, u_h, f_h)$ can be applied to the fine grid system (\ref{fe_equation})
in the solution of which we are interested. The MG iteration with procedure MGM is repeated on the fine grid until the fine residual
 is considered small enough, see Algorithm \ref{alg:fmg} for details.

\begin{figure}
  \centering
  \includegraphics[width=5in]{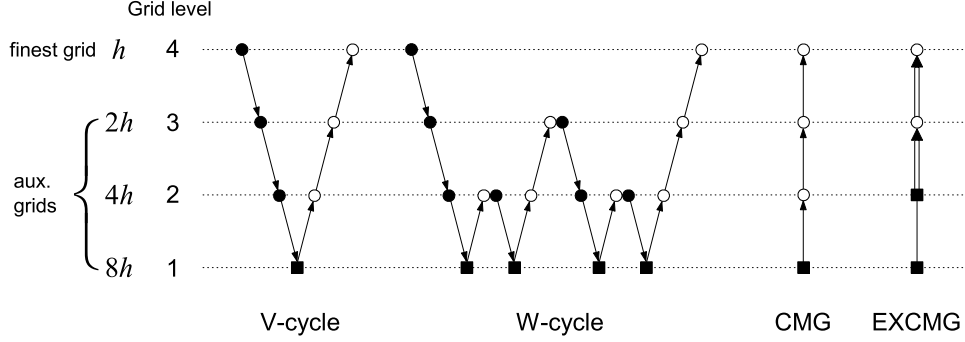}\\
  \caption{The four level structure of the V- and W-cycles, CMG and ECMG methods. In the diagram, $\bullet$ denotes pre-smoothing steps,
  $\circ$ denotes post-smoothing steps, $\uparrow$ denotes prolongation (usually defined by linear interpolation), $\downarrow$ denotes restriction,
  $\Uparrow$ denotes extrapolation and quadratic interpolation, and $\blacksquare$ denotes direct solver.
  }\label{VW}\end{figure}

To fully specify the Algorithm \ref{alg:mg}, one needs to specify the numbers
of pre- and post-smoothing steps on each grid $N_{pre}$ and  $N_{post}$,
the restriction (fine-to-coarse) operator $R_h^{2h}$, the prolongation (coarse-to-fine) operation $P_{2h}^h$, and the recursion parameter $N_{cycles}$.
One iteration of a MG method, from the finest grid to the coarsest grid and back
to the finest grid again, is called a cycle. The exact structure of MG cycle depends on the parameter $N_{cycles}$,
which is the number of two-grid iterations at each intermediate stage, greater or equal to one. With $N_{cycles}=1$,
the so-called V-cycle is generated, while $N_{cycles}=2$ leads to W-cycle. The four level structures of the V- and W-cycles are
illustrated in Fig.\ref{VW}

The proper choices of both residual restriction and error prolongation operators play a key role in the successful development of the MG
method.
For prolongation operator $P_{2h}^h$, we use the standard tri-linear interpolation operators which preserves the symmetry of the true solution exactly. And the 27-point fully weighted residual restriction
is used in this study, which is a natural extension from the 2D 9-point full weighting restriction to 3D space. The restriction
operator $R_h^{2h}$ is defined as follows \cite{Trottenberg2001}:
\begin{equation}\label{3Dres}
  R_h^{2h} u_h(x,y,z) = \sum_{i,j,k=-1}^{+1} w_{ijk}u_h(x+i\Delta x, y+j\Delta y, z+k\Delta z),
\end{equation}
where
\begin{equation}
  w_{ijk} = \left\{  \begin{aligned}
         &1/8          &\textrm{if} \quad &i=j=k=0,\\
         &1/16         &\textrm{if} \quad &|i|+|j|+|k|=1,\\
         &1/32         &\textrm{if} \quad &|i|+|j|+|k|=2,\\
         &1/64         &\textrm{if} \quad &|i|+|j|+|k|=3.\\
         \end{aligned} \right.
\end{equation}
In this manner, the weighting satisfied the conservation property of the integrals.

\section{Extrapolation cascadic multigrid methods}

The CMG method proposed by Deuflhard and Bornemann in \cite{Bornemann1996} is a type of method
which requires no coarse grid correction at all that may be viewed as a one-way MG (see Figure \ref{VW}).
Because no coarse grid correction is needed, computational cost is greatly reduced compared to the classical MG methods.
Since the 1990s, the method received quite a bit of attention from researchers because of its high efficiency and simplicity \cite{Shaidurov1996, Braess1999,Timmermann2000,Shaidurov2000a,Shaidurov2000,Shi2000,Shi2001,Braess2002,Stevenson2002,Bi2004,Shi2004,Zhou2004,Man2006,Shi2007,Du2008,Wang2008,Xu2009,Yu2011}.
By using Richardson extrapolation
and bi-quadratic interpolation techniques, an extrapolation cascadic multigrid (ECMG) method for elliptic problems was first proposed
in 2008 by Chen et al. in \cite{Chen2008,Hu2009}.
And some 2D numerical experiments are presented in \cite{Chen2011} to show the advantage of ECMG in comparison with CMG method.

\subsection{Description of the ECMG algorithm}

The key ingredients of the ECMG method are extrapolation and quadratic interpolation (see Figure \ref{VW}), which are used
to provide a better initial guess of the iterative solution on the next finer grid than one obtained by using linear interpolation in CMG. The ECMG
algorithm is given as follows \cite{Chen2011}:

\begin{algorithm}[htb]         

\caption{ ECMG : $u_h$ $\Leftarrow$ ECMG$(A_h, f_h, L, m_j)$}             

\label{alg:oECMG}                  

\begin{algorithmic}[1]                

\STATE $u_H$ $\Leftarrow$ DSOLVE($A_H u_H=f_H$)
\STATE $u_{H/2}$ $\Leftarrow$ DSOLVE($A_{H/2} u_{H/2}=f_{H/2}$)
\STATE $h=H/2$

  \FOR {$i=1$ to $L$}
   \STATE $h = h/2$
   \STATE ${u}_{h} = \textrm{EXP}_{finite}(u_{2h}, u_{4h})$
    \FOR {$j=1$ to $m_j$}
        \STATE $u_h \Leftarrow$ CG$(A_h, u_h, f_h)$
    \ENDFOR
  \ENDFOR
\end{algorithmic}

\end{algorithm}

In Algorithm \ref{alg:oECMG}, procedure CG denotes conjugate gradient iteration, and procedure EXP$_{finite}$ denotes extrapolation and bi-quadratic interpolation operator on three levels of the embedded rectangular grids,
which is used to get a good initial guess of the iterative solution on the next finer grid, see line 6 and line 8 in the Algorithm \ref{alg:oECMG}.
The positive integer $L$, total number of grids except first two embedded grids, indicates that
the mesh size of the finest grid is $\frac{H}{2^{L+1}}$. And on the first two coarse grids, a direct solver DSOLVE is used since the size of the linear system is small, see line 1-2 in the Algorithm \ref{alg:oECMG}.

The number of iterations on each grid is defined by
\begin{equation}\label{ee}
  m_j=\lceil m_L\cdot \beta^{L-j}\rceil,\quad j=2,3,\cdots,L,
\end{equation}
where $m_L$ is the iteration count on the finest grid, $\beta$ is usually set to be between 2 and $2^d$,
 and $\lceil x \rceil$ denotes the smallest integer greater than or equal to $x$. In \cite{Chen2009}, authors compared different settings for $m_L$ and $\beta$, and found that ECMG with $m_j = 4\times 4^{L-j}$  can achieve high accuracy with quick convergence.
In 2007, Shi, Xu and Huang \cite{Shi2007} proposed an economical CMG method,
which uses a new criteria for choosing the smoothing steps on each level, when $d = 2$, that is
\begin{description}
  \item[(\rmnum{1})] if $j>L_0$, then $ m_j=\lceil m_L\cdot \beta^{L-j}\rceil$,
  \item[(\rmnum{2})] if $j\leq L_0$, then $m_j=\lceil (L-(2-\epsilon_0)j)h_j^{-2} \rceil$.
\end{description}
Here $0<\epsilon_0\leq 1$ is a fixed positive number, and $m_L = m_0(L-L_0)^2$. Note that in standard cascadic multigird algorithm
$m_L=m_0L^2$. The level parameter $L_0$, which depends on $L, m_L, \beta$ and the size of the coarest level $H$,
can be defined as the largest integer satisfying
\begin{equation}
  L_0\leq \min\left\{\frac{L\log\beta+\log m_L+4\log H}{\log \beta+4\log 2},\frac{L}{2}\right\}.
\end{equation}

To our best knowledge, the ECMG method has mainly been applied to solve 2D elliptic boundary value problems, but it is very important to solve the 3D elliptic
 problems efficiently and accurately in many engineering areas and the extension of ECMG method from 2D to 3D is nontrivial.
The advantage of the ECMG method will be much more reflected when it is applied to 3D problems since the size of linear discrete system is much larger.
Next, we propose a new ECMG method for solving 3D elliptic problems which is stated in the following algorithm~\ref{alg:ECMG}:

\begin{algorithm}[!htb]         

\caption{New ECMG : $(u_h, \tilde{u}_h)$ $\Leftarrow$ ECMG($A_h, f_h, L ,\epsilon$)}             

\label{alg:ECMG}                  

\begin{algorithmic}[1]                

\STATE $u_H$ $\Leftarrow$ DSOLVE($A_H u_H=f_H$)
\STATE $u_{H/2}$ $\Leftarrow$ DSOLVE($A_{H/2} u_{H/2}=f_{H/2}$)
\STATE $h=H/2$

  \FOR {$i=1$ to $L$}
   \STATE $h = h/2$
   \STATE ${w}_{h} = \textrm{EXP}_{finite}(u_{2h}, u_{4h})$  $\quad \quad \quad\rhd$ $u_h=w_h$ is used as the initial guess  for JCG solver
    \WHILE {$||A_h u_h -f_h||_2>\epsilon \cdot ||f_h||_2 $}
        \STATE $u_h \Leftarrow$ JCG$(A_h, u_h, f_h)$
    \ENDWHILE
    \STATE $\tilde{u}_{h} = \textrm{EXP}_{true}(u_{h}, u_{2h})$ $\quad \quad \quad \quad \rhd$ $\tilde{u}_{h}$ is a higher-order approximate solution
  \ENDFOR
\end{algorithmic}
\end{algorithm}

In Algorithm \ref{alg:ECMG}, $\textrm{EXP}_{finite}(u_{2h}, u_{4h})$ denotes a third order approximation to the FE solution $u_h$ obtained by Richardson extrapolation and tri-quadratic Serendipity interpolation techniques from the numerical solutions on embedded hexahedral grids, while $\textrm{EXP}_{true}(u_{h}, u_{2h})$
denotes a fourth order extrapolated solution on fine grid with mesh size $h$ from two second order numerical solutions $u_h$ and $u_{2h}$.
The detail procedures of extrapolation and Serendipity interpolation on  embedded hexahedral grids
are described in the next Section \ref{extra}. The differences between our new ECMG and the existing ECMG are illustrated as follows:
\begin{enumerate}[(1)]
  \item Instead of performing a fixed number of CG iterations (see line 7 in Algorithm \ref{alg:oECMG})
  as used in the usual CMG methods, a tolerance $\epsilon$  related to the relative residual is used in our iterative solvers (see line 7 in Algorithm \ref{alg:ECMG}),
  which enables us to not only avoid the difficulty of choosing the number of iterations on each grid $m_i$,
  but also allows us to conveniently obtain the numerical solution with the desired accuracy.

 \item Algorithm \ref{alg:ECMG} takes JCG as the MG smoother, which is written as the line 8 in the Algorithm \ref{alg:ECMG}.
 The Jacobi preconditioner, which consists only of the diagonal elements of the matrix,
can be constructed from the coefficient matrix without any more extra work. Although JCG requires a little additional work than CG,
JCG smoother typically yields better convergence properties
than CG smoother (see Table~\ref{iterations8}-~\ref{iterations10} in Section \ref{sec-nu}).

  \item A fourth order extrapolated solution $\tilde{u}_{h}$ (see Table~\ref{table1}-~\ref{table4} for details) on the current grid based on numerical solutions of two level grids, current grid and previous (coarser) grid, is constructed in Algorithm \ref{alg:ECMG} (see line 10)
   in order to enhance the accuracy of the numerical solution $u_h$.
\end{enumerate}

In the following, we denote the above new  ECMG method as $\textrm{ECMG}_{jcg}$ and denote $\textrm{ECMG}_{cg}$ for the Algorithm \ref{alg:ECMG} where the line 8 is replaced by $u_h \Leftarrow$ CG$(A_h, u_h, f_h)$.

\subsection{Extrapolation and quadratic interpolation}\label{extra}
It is well known that the extrapolation method, which was established by Richardson
in 1926, is an efficient procedure for increasing the solution accuracy of many
problems in numerical analysis. In 1983, Marchuk and Shaidurov \cite{Marchuk1983}
studied systematically the application of this method in finite difference method.
Since then, this technique has been well demonstrated in the frame of the
FE method \cite{Lin1983,Lin1985,Blum1986,Chen1989,Fairweather2006,Asadzadeh2010}.

In this Section, we will explain how to use extrapolation and interpolation techniques to obtain a high order accurate approximate solution to the problem (\ref{bvp}), and a high order accurate approximation to the FE solution. And it can be regarded as another important application of the extrapolation method to provide a good initial guess of iterative solution.

\subsubsection{Extrapolation for the true solution}
For simplicity, let us first consider the three levels of embedded grids $Z_i(i=0,1,2)$ with mesh-sizes $h_i=h_0/2^i$ in one dimension.
Let $e^i=U^i-u$ be the error of the linear FE solution $U^i\in V_{h_i}$,
it can be proved that under certain smoothness conditions the error at node $x_k$
has the the following asymptotic expansion
\begin{equation}\label{asy}
  e^i(x_k)=A(x_k)h_i^2+O(h_i^4),
\end{equation}
where $A(x)$ is a suitably smooth function independent of $h_i$.

It is well known that the extrapolation methods can only offer a fourth-order accurate approximation to the true solution
at coarse grid points. The Richardson extrapolation formula at coarse grid points $\{x_j,x_{j+1}\}$ can be written as
\begin{equation}\label{t1}
  \widetilde{U}_k^1 := \frac{4 U^1_k - U^0_k}{3} = u(x_k) + O(h_0^4),\ \ k=j,j+1.
\end{equation}

In fact, by using the linear interpolation formula, one can also obtain a fourth-order accurate approximation at fine grid points.
Chen and Lin \cite{Chen1989} proposed a fine grid extrapolation formula
\begin{equation}\label{chenlin}
  \widetilde{U}_{j+1/2}^1 := U_{j+1/2}^1 + \frac{1}{6}(U^1_j - U^0_j + U^1_{j+1} - U^0_{j+1}) = u(x_{j+1/2}) + O(h_0^4),
\end{equation}
which can give directly the higher order accuracy at fine grid points.

From (\ref{asy}), we can obtain
\begin{equation}\label{aaa}
  A(x_k) = \frac{4}{3h_0^2}(U_k^0-U_k^1) + O(h_0^2),\ \ k=j,j+1.
\end{equation}
Applying the error estimate of the linear interpolation,
\begin{equation}\label{A_12}
  A(x_{j+1/2})=\frac{1}{2}(A(x_j)+A(x_{j+1}))+O(h_0^2).
\end{equation}
Substituting eq. (\ref{aaa}) into eq. (\ref{A_12}) yields
\begin{equation}\label{Axmid}
  A(x_{j+1/2}) = \frac{2}{3h_0^2}(U_j^0-U_j^1) + \frac{2}{3h_0^2}(U_{j+1}^0-U_{j+1}^1) + O(h_0^2).
\end{equation}
Since
\begin{equation}
  U_{j+1/2}^1 = u(x_{j+1/2}) + \frac{1}{4} A(x_{j+1/2})h_0^2 + O(h_0^4),
\end{equation}
by using (\ref{Axmid}), the extrapolation formula (\ref{chenlin}) is obtained.

\subsubsection{Extrapolation for the FE solution}
Next, we will explain,
given the FE solutions $U^0$ and $U^1$, how to use the extrapolation method to obtain a third order (to be proved in subsection \ref{error}) approximation $W^2$ for the { FE solution $U^2$} rather than the { true solution $u$}.

\begin{figure}
   \centering
   \scalebox{0.5}{\includegraphics{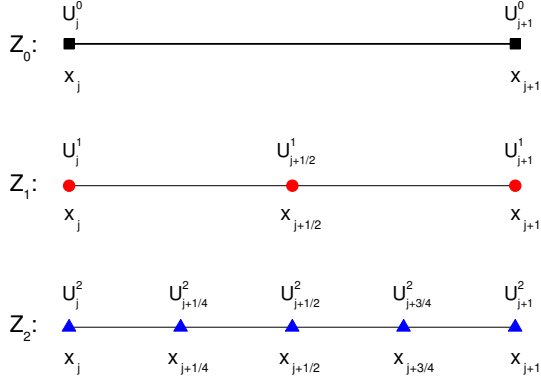}}
   \caption{Three embedded grids in 1D.}\label{Fig:1}
\end{figure}

Adding one midpoint and two four equal division points,
the coarse mesh element $(x_j,x_{j+1})$ is uniformly refined into four elements as shown in Fig. \ref{Fig:1}.
And we obtain a five points set $$\Big\{x_j,x_{j+1/4},x_{j+1/2},x_{j+3/4},x_{j+1}\Big\}$$
belonging to fine mesh $Z_2$.


Assume there exists a constant $c$ such that
\begin{equation}\label{comb}
  c U^1+(1-c) U^0=U^2+O(h_0^4),
\end{equation}
that is to say, we use a linear combination of $U^0$ and $U^1$
to approximate {the FE solution $U^2$.}
Substituting the asymptotic error expansion (\ref{asy}) into (\ref{comb}), we obtain $c=5/4$ and an
extrapolation formula
\begin{equation}\label{jiedian}
   W^2_k := \frac{5 U^1_k - U^0_k}{4} = U^2_k+O(h_0^4),\ \ k=j,j+1,
\end{equation}
 at nodes $x_j$ and $x_{j+1}$. To derive the extrapolation formula at mid-point $x_{j+1/2}$,  we
have
\begin{equation}\label{bbb}
  U_{j+1/2}^2 = U_{j+1/2}^1 - \frac{3}{16} A(x_{j+1/2}) h_0^2 + O(h_0^4),
\end{equation}
by using equation (\ref{asy}).
Substituting eq.(\ref{Axmid}) into eq.(\ref{bbb}) yields the following mid-point extrapolation formula,
\begin{align}
      W^2_{j+1/2}:=U^1_{j+1/2}+ \frac{1}{8}(U^1_j - U^0_j + U^1_{j+1} - U^0_{j+1})={U}^2_{j+1/2} + O(h_0^4). \label{zhongdian}
\end{align}
Once the three initial values $W_j^2,W_{j+1}^2$ and $W_{j+1/2}^2$  are obtained, we can get the following four equal division point extrapolation formulas by using
quadratic interpolation,
\begin{align}
      W^2_{j+1/4}&:=\displaystyle\frac1{16}\big[(9U^1_{j}+12U^1_{j+1/2}-U^1_{j+1})-(3U^0_j+U^0_{j+1})\big],\label{sifen1}\\
      W^2_{j+3/4}&:=\displaystyle\frac1{16}\big[(9U^1_{j+1}+12U^1_{j+1/2}-U^1_j)-(3U^0_{j+1}+U^0_j)\big].\label{sifen2}
\end{align}

 \subsection{Three dimensional case}

\begin{figure}
  \centering
  \includegraphics[width=1.0\textwidth]{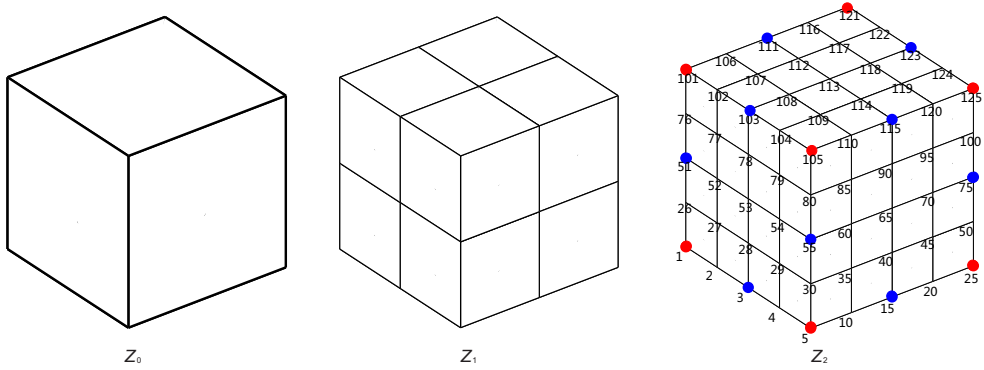}
  \caption{ Three embedded grids in 3D.
  }\label{Fig3D}
  \end{figure}
In this subsection, we explain how to obtain a third order accurate
 approximation ${W}^2$ to the FE solution $U^2$, and a fourth order accurate approximate solution $\tilde{U}^1$ to the problem (\ref{bvp})
 for embedded hexahedral grids as shown in Fig. \ref{Fig3D}.

The construction process of the approximation $W^2$ are as follows:
\begin{description}

  \item[Corner Nodes (1, 5, 21, 25, 101, 105, 121, 125):]
   The approximate values at 8 corner nodes `\textcolor[rgb]{1,0,0}{$\medbullet$}' can be obtained by using the
extrapolation formulae (\ref{jiedian})

  \item[Midpoints in $x$-direction (3, 23, 103, 123):]
    The approximate values at these 4 midpoints `\textcolor[rgb]{0,0,1}{$\medbullet$}' can be obtained by using the mid-point extrapolation formulae (\ref{zhongdian}) in $x$-direction.

    \item[Midpoints in $y$-direction (11, 15, 111, 115):]
    The approximate values at these 4 midpoints `\textcolor[rgb]{0,0,1}{$\medbullet$}' can be obtained by using the mid-point extrapolation formulae (\ref{zhongdian}) in $y$-direction.

    \item[Midpoints in $z$-direction (51, 55, 71, 75):]
    The approximate values at these 4 midpoints `\textcolor[rgb]{0,0,1}{$\medbullet$}' can be obtained by using the mid-point extrapolation formulae (\ref{zhongdian}) in $z$-direction.

\item[Other 105 points:] The approximate values of remaining 105 $(5^3-20)$ grid points can be obtained by using tri-quadratic Serendipity interpolation with the known 20-node (8 corner nodes and 12 mid-side nodes) values.
\end{description}

The tri-quadratic Serendipity interpolation function in terms of natural coordinates ($\xi,\eta,\zeta$) is
    \begin{equation}\label{inter}
      W(\xi,\eta,\zeta)=\sum_i N_i(\xi,\eta,\zeta) W_i,
    \end{equation}
    where the shape functions $N_i$ can be written as follows \cite{Liu2013}:
  \begin{equation}
    {N_i}(\xi ,\eta ,\zeta ) = \left\{ \begin{aligned}
  &\frac{1}{8}(1 + {\xi_i\xi})(1 + {\eta_i\eta})(1 + {\zeta_i\zeta})({\xi_i\xi} + {\eta_i\eta} + {\zeta_i\zeta} - 2),\quad &i& = 1,5,21,25,101,105,121,125  \\
  &\frac{1}{4}(1 - {\xi ^2})(1 + {\eta_i\eta})(1 + {\zeta_i\zeta}),&i& = 3,23,103,123  \\
  &\frac{1}{4}(1 - {\eta ^2})(1 + {\xi_i\xi})(1 + {\zeta_i\zeta}),&i &= 11,15,111,115  \\
  &\frac{1}{4}(1 - {\zeta ^2})(1 + {\xi_i\xi})(1 + {\eta_i\eta}),&i& = 51,55,71,75  \\
\end{aligned}  \right.
  \end{equation}
  where $(\xi_i,\eta_i,\zeta_i)$ is the natural coordinate of node $i$. For example,

 \begin{equation}
     {W_{35}} = \sum\limits_i^{} {{N_i}(1, - 0.5, - 0.5){W_i}}  = \frac{9}{{16}}{W_{15}} - \frac{3}{{16}}{W_{25}} + \frac{9}{{16}}{W_{55}} + \frac{3}{{16}}{W_{75}} - \frac{3}{{16}}{W_{105}} + \frac{3}{{16}}{W_{115}} - \frac{1}{8}{W_{125}}.
   \end{equation}

\begin{remark}
  Extrapolation and tri-quadratic interpolation defined by Eq.(\ref{jiedian}), Eq.(\ref{zhongdian}) and Eq.(\ref{inter})
  are local operations, which can be be done very cheaply and effectively.
\end{remark}

\begin{remark}
Since the center of each face such as the point labeled 53 is the midpoint of two face diagonals,
we can also use the mid-point extrapolation formulae (\ref{zhongdian}) to obtain two different extrapolation values,
 and take their arithmetic mean as the approximated value. A similar procedure can be done for the center point of the cube labeled 63,
 which can be viewed as the midpoint of space diagonals.
 Then the approximated values of remaining 98 $(5^3-27)$ points can be obtained
 by using tri-quadratic Lagrange interpolation with the known 27-node values. This tri-quadratic Lagrange interpolation
 also gives a third order approximation of the FE solution on the finer grid followed
 by a similar error analysis in the next subsection, in addition,
 numerical results show that it has almost the same accuracy as the tri-quadratic Serendipity interpolation with the known 20 nodes.
  Thus, we only present the numerical results by using the method of tri-quadratic Serendipity interpolation in the next Section.
\end{remark}

When constructing the fourth order accurate approximate solution  $\tilde{U}^1$ based on two second order accurate solutions $U^0$ and $U^1$,
the Richardson extrapolation formula (\ref{t1}) can be directly used for 8 coarse grid points,
while the fourth order extrapolation formula (\ref{chenlin}) can be used for the 12 centers of  edges,
the 6 centers of faces as well as the center point of the coarse hexahedral element.

\begin{remark}
In the existing literature,  there are already some extrapolation related methods which use the two numerical solutions from
the two level of grids to obtain high order accuracy solution on the fine grid. For example,
Wang et al. \cite{Wang2009} first used the Richardson extrapolation for two  fourth order compact finite difference solutions on
two level of grids to achieve sixth order accuracy on the even number of grids of the fine grid, however,  the method of~\cite{Wang2009}
relies on an operator based interpolation iterative strategy for odd number of grids in order to achieve the sixth order accuracy on the whole fine grid,
which increases the computational cost. Our extrapolation formulas (\ref{t1})-(\ref{chenlin}) can be used to obtain fourth order accurate solution on the
whole fine grid from two second order accurate solutions on two level of grids directly and  cheaply without any extra work.
\end{remark}

\subsection{The error analysis of initial guess $W^2$}\label{error}
Let $e=W^2-{U}^2$ be the approximation error of extrapolation formulas, and
assume that it has a continuous derivative up to order 3 on interval $[x_j,x_{j+1}]$.
From (\ref{jiedian}) and (\ref{zhongdian}) we obtain the equation
 \begin{equation}\label{err}
   e(x_k) = O(h_0^4),\quad k=j,j+1/2,j+1.
 \end{equation}
From polynomial interpolation theory, the
error of quadratic interpolation $I_2f$ can be represented as,
\begin{equation}
  R_2(x) = e - I_2e = \frac{1}{3!}(x-x_j)(x-x_{j+1/2})(x-x_{j+1})e^{(3)}(\xi),
\end{equation}
where $\xi$ is a point of $[x_j,x_{j+1}]$ that depends on $x$.
Especially at four equal division points we have
\begin{equation}\label{si1}
  R_2(x_{j+1/4}) = \frac{h_0^3}{128}e^{(3)}(\xi_1) = \frac{h_0^3}{128}e^{(3)}(x_{j+1/2}) + O(h_0^4),
\end{equation}
and
\begin{equation}\label{si2}
  R_2(x_{j+3/4}) = -\frac{h_0^3}{128}e^{(3)}(\xi_2) = -\frac{h_0^3}{128}e^{(3)}(x_{j+1/2}) + O(h_0^4)\approx - R_2(x_{j+1/4}).
\end{equation}
It follows from eq. (\ref{si1}), eq.(\ref{si2}) and eq.(\ref{err}) that
\begin{equation}\label{err_si}
   e(x_k) = I_2 e(x_k) + R_2(x_k) = O(h_0^3) ,\quad k=j+1/4,j+3/4,
 \end{equation}
 which means that four equal division point extrapolation formulas (\ref{sifen1}) and (\ref{sifen2}) are only  third-order approximations.

The above error analysis can be directly extended to 3D case (also see numerical verification of Table \ref{table11}-\ref{table5} in the next Section).
In addition, eq.(\ref{si2}) implies that the  error $e(x)$ forms a high-frequency oscillation in whole domain,
however, it can be smoothed out after very few iterations (see Fig. \ref{Fig.lable} for details).

\section{Numerical experiments}\label{sec-nu}
%

In this Section, we will illustrate the  efficiency  of $\textrm{ECMG}_{jcg}$ by comparing to  $\textrm{ECMG}_{cg}$
as well as the classical V-cycle and W-cycle MG methods,
and present numerical results for three examples with smooth and singular solutions obtained by the most efficient algorithm $\textrm{ECMG}_{jcg}$.
In our experiments, we use the most popular smoothing method Gauss-Seidel smoother as the relaxation smoother of classical MG methods,
which usually leads to a good convergence rate.
Our code is written in Fortran 90 and compiled with  Intel Visual Fortran Compiler XE 12.1 compiler. All programs are carried out on a server with Intel(R) Xeon(R) CPU E5-2680  (2.80GHz)  and 32G RAM.

\begin{table}
\tabcolsep=15pt
\caption{Comparison of the numbers of iterations and CPU times (in seconds) of four different MG methods with $\epsilon = 10^{-8}$ for Problem 1.} \centering
\begin{threeparttable}
\begin{tabular}{ccccc}
    \hline
Mesh & ECMG$_{jcg}$  &   ECMG$_{cg}$  &    V(1,1)\tnote{1}     &   W(2,1)\tnote{2}     \\
\hline
$  32\times  32\times  32$ &  7 &   58&  26  &  432   \\
$  64\times  64\times  64$ & 10 &   82&  26  &  216   \\
$ 128\times 128\times 128$ & 18 &   93&  26  &  108   \\
$ 256\times 256\times 256$ &  3 &   58&  26  &  54   \\
$ 512\times 512\times 512$ &  3 &    9&  26  &  27   \\
\hline
CPU (s) &  {147}  &  {169}  & {452}  &  {466}\\
\hline
\end{tabular}
      \begin{tablenotes}
        \footnotesize
        \item[1] The number of MG V(1,1) cycle is 13.
        \item[2] The number of MG W(2,1) cycle is 9.
      \end{tablenotes}
    \end{threeparttable}
\label{iterations8}
\end{table}

\begin{table}
\tabcolsep=15pt
\caption{Comparison of the numbers of iterations and CPU times (in seconds) of four different MG methods with $\epsilon = 10^{-9}$ for Problem 1.} \centering
\begin{threeparttable}
\begin{tabular}{ccccc}
    \hline
Mesh & ECMG$_{jcg}$  &   ECMG$_{cg}$  &    V(1,1)\tnote{1}     &   W(2,1)\tnote{2}     \\
\hline
$  32\times  32\times  32$ &  8 &   71&  30  & 480   \\
$  64\times  64\times  64$ &  9 &   98&  30  & 240   \\
$ 128\times 128\times 128$ & 16 &  150&  30  &  120   \\
$ 256\times 256\times 256$ & 78 &  162&  30  &  60   \\
$ 512\times 512\times 512$ &  3 &   50&  30  &  30   \\
\hline
CPU (s) &  {162}  &  {284}  & {489}  &  {545}\\
\hline
\end{tabular}
      \begin{tablenotes}
        \footnotesize
        \item[1] The number of MG V(1,1) cycle is 15.
        \item[2] The number of MG W(2,1) cycle is 10.
      \end{tablenotes}
    \end{threeparttable}
\label{iterations9}
\end{table}

{\bf Problem 1.} The test Problem 1 can be written as
\begin{equation}\label{prob1}
         \frac{\partial^2 u}{\partial x^2}+\frac{\partial^2 u}{\partial y^2}+\frac{\partial^2 u}{\partial z^2} =-\frac{3}{4}\pi^2 \sin(\frac{\pi}{2} x)\sin(\frac{\pi}{2} y)\sin(\frac{\pi}{2} z),\quad  \textrm{in } \Omega=[0,1]^3,
\end{equation}
where the boundary conditions are
\begin{equation}
  u(0,y,z)=u(x,0,z)=u(x,y,0)=0,\quad \frac{\partial u}{\partial n}(1,y,z)=\frac{\partial u}{\partial n}(x,1,z)=\frac{\partial u}{\partial n}(x,y,1)=0.
\end{equation}
The analytic solution of eq. (\ref{prob1}) is
$$u(x,y,z)=\sin(\frac{\pi}{2} x)\sin(\frac{\pi}{2} y)\sin(\frac{\pi}{2} z),$$
where the exact solution is a sufficiently smooth function which has arbitrary order smooth derivatives and the variation of the function is the same in three directions.

\begin{table}
\tabcolsep=15pt
\caption{Comparison of the numbers of iterations and CPU times (in seconds) of four different MG methods with $\epsilon = 10^{-10}$ for Problem 1.} \centering
\begin{threeparttable}
\begin{tabular}{ccccc}
    \hline
Mesh & ECMG$_{jcg}$  &   ECMG$_{cg}$  &    V(1,1)\tnote{1}     &   W(2,1)\tnote{2}     \\
\hline
$  32\times  32\times  32$ &   9 &   76&  32  &  528   \\
$  64\times  64\times  64$ &   9 &  124&  32  &  264   \\
$ 128\times 128\times 128$ &  26 &  172&  32  &  132   \\
$ 256\times 256\times 256$ & 112 &  258&  32  &  66   \\
$ 512\times 512\times 512$ &   3 &   80&  32  &  33    \\
\hline
CPU (s) &  {177}  &  {371}  & {522}  &  {534}\\
\hline
\end{tabular}
      \begin{tablenotes}
        \footnotesize
        \item[1] The number of MG V(1,1) cycle is 16.
        \item[2] The number of MG W(2,1) cycle is 11.
      \end{tablenotes}
    \end{threeparttable}
\label{iterations10}
\end{table}

Using 7 embedded grids with the coarsest grid $8\times8\times8$, we give the CPU time and the number of iterations at each level of grid for $\textrm{ECMG}_{jcg}$, $\textrm{ECMG}_{cg}$, the classical V-cycle and W-cycle MG methods using three different tolerances $\epsilon=10^{-8}, 10^{-9}$ and $10^{-10}$  in  Table~\ref{iterations8},  Table~\ref{iterations9} and  Table~\ref{iterations10}, respectively. Here the CPU time is total computational time which not only includes the time to solve the linear system by MG methods but also includes all other necessary computational time in the program. Since the first two coarse level of grids are solved by direct solver, we only list the results start from the third level of grid $32\times32\times32$.

As we can see from Table~\ref{iterations8} to Table~\ref{iterations10} that the  $\textrm{ECMG}_{jcg}$ uses the minimum number of iterations and minimum CPU time and has the best efficiency  among all MG methods,  $\textrm{ECMG}_{cg}$ is a little worse than $\textrm{ECMG}_{jcg}$, while the classical MG with W-cycle has the worst efficiency and has the similar efficiency as the classical MG with V-cycle. In addition, we can see that when the control parameter $\epsilon$ becomes smaller, the number of iterations and the CPU time for all MG methods increase as expected, and the efficiency of $\textrm{ECMG}_{jcg}$ becomes much better than $\textrm{ECMG}_{cg}$ while the classical MG with V-cycle still has the similar efficiency as the classical MG with W-cycle.
It should be noted that when the tolerance becomes $10^{-10}$, only 3 JCG iterations are required on the finest grid $512\times 512 \times 512$ (more than 135 million unknowns) for ECMG$_{jcg}$ method, which costs less than 3 minutes (see second column of Table \ref{iterations10}).
Moreover, numerical comparisons for efficiency are also carried out for the next two examples (not listed in the paper), results show that $\textrm{ECMG}_{jcg}$ always has the best efficiency among these four MG methods.

For the sake of clarity, however, it should be mentioned that the process of solving the linear system only  required a small portion of the total computing time compared to the effort spent in discretizing the problem, i.e., computing the stiffness matrices and load vectors.
Additionally, since we are using different MG methods to solve the same linear system, the accuracy of all methods should be similar under the same tolerance. Thus, in the following part of this paper, we present only the numerical results obtained by $\textrm{ECMG}_{jcg}$ since it is the most efficient method.


\begin{table}[t]
\tabcolsep=6pt
\caption{Errors and convergence rates with $\epsilon=10^{-8}$ in $L_{2}$ norm for Problem 1.} \centering
\begin{tabular}{|c|c|c|cc|cc|cc|c|}
    \hline
    \multirow{2}{*}{Mesh}& \multirow{2}{*}{Iters} & \multirow{2}{*}{RRe} & \multicolumn{2}{c|}{$||U_h-u||_{2}$} & \multicolumn{2}{c|}{$||\widetilde{U}_h-u||_{2}$} & \multicolumn{2}{c|}{$||W_h-U_h||_{2}$}  &  \multirow{2}{*}{$r_h$} \\
    \cline{4-9}&   &   & Error &   Order &   Error &   Order  &   Error &   Order & \\
\hline
$  32\times  32\times  32$ &   7 & $9.99(-9)$ &$1.42(-4)$ &      &  $1.96(-7) $  &       &  $2.54(-5)$ &      &   $1.79(-1)$   \\
$  64\times  64\times  64$ &  10 & $8.14(-9)$ &$3.55(-5)$ & 2.00 &  $1.24(-8) $  & 3.98  &  $3.18(-6)$ & 2.99 &   $8.96(-2)$   \\
$ 128\times 128\times 128$ &  18 & $7.68(-9)$ &$8.87(-6)$ & 2.00 &  $7.83(-10)$  & 3.99  &  $3.99(-7)$ & 3.00 &   $4.50(-2)$   \\
$ 256\times 256\times 256$ &   3 & $6.26(-9)$ &$2.22(-6)$ & 2.00 &  $4.74(-10)$  & 0.72  &  $4.99(-8)$ & 3.00 &   $2.25(-2)$   \\
$ 512\times 512\times 512$ &   3 & $6.00(-9)$ &$5.55(-7)$ & 2.00 &  $4.69(-10)$  & 0.02  &  $6.25(-9)$ & 3.00 &   $1.13(-2)$   \\
\hline
\end{tabular}\label{table11}
\end{table}

\begin{table}[t]
\tabcolsep=6pt
\caption{Errors and convergence rates with $\epsilon=10^{-8}$ in $L_{\infty}$ norm for Problem 1.} \centering
\begin{tabular}{|c|c|c|cc|cc|cc|}
    \hline
    \multirow{2}{*}{Mesh}& \multirow{2}{*}{Iters} & \multirow{2}{*}{RRe} & \multicolumn{2}{c|}{$||U_h-u||_{\infty}$} & \multicolumn{2}{c|}{$||\widetilde{U}_h-u||_{\infty}$} & \multicolumn{2}{c|}{$||W_h-U_h||_{\infty}$}  \\
    \cline{4-9}&   &   & Error &   Order &   Error &   Order  &   Error &   Order  \\
\hline
$  32\times  32\times  32$ &   7 & $9.99(-9)$ &$4.02(-4)$ &      &  $1.11(-6) $ &       &  $2.54(-5)$ &        \\
$  64\times  64\times  64$ &  10 & $8.14(-9)$ &$1.00(-4)$ & 2.00 &  $6.95(-8) $ & 4.00  &  $3.18(-6)$ & 2.99   \\
$ 128\times 128\times 128$ &  18 & $7.68(-9)$ &$2.51(-5)$ & 2.00 &  $4.39(-9)$  & 3.98  &  $3.99(-7)$ & 3.00   \\
$ 256\times 256\times 256$ &   3 & $6.26(-9)$ &$6.28(-6)$ & 2.00 &  $1.56(-9)$  & 1.50  &  $4.99(-8)$ & 3.00   \\
$ 512\times 512\times 512$ &   3 & $6.00(-9)$ &$1.57(-6)$ & 2.00 &  $1.38(-9)$  & 0.17  &  $6.25(-9)$ & 3.00   \\
\hline
\end{tabular}\label{table22}
\end{table}

We present the numerical results for Problem 1  obtained by $\textrm{ECMG}_{jcg}$ with $\epsilon=10^{-8}$ in Table~\ref{table11}-\ref{table22},
and with $\epsilon=10^{-9}$ in Table~\ref{table1}-\ref{table2}, where ``Iters'' denotes the number of iterations needed for the JCG solver to achieve that  the relative residual is less than the given tolerance $\epsilon$ while ``RRe'' denotes the corresponding relative residual. And the same notations are used in all the following tables.
Table~\ref{table11} and  Table~\ref{table1} list the number of iterations, the relative residual of the numerical solution on each grid, the $L_2$ error between the FE solution $U_h$ and the exact solution $u$, the $L_2$ error between the extrapolated solution $\tilde{U}_h$ and the exact solution $u$, the $L_2$ error between the initial  guess $W_h$ and the FE solution $U_h$, all convergence rates, and  the ratio $r_h$ defined as
\begin{equation}
r_h=\frac{\parallel W_h-U_h\parallel_{2}}{\parallel U_h-u\parallel_{2}},
\end{equation}
which is used to measure how good $W_h$ approximates $U_h$.  While Table~\ref{table22} and  Table~\ref{table2} give all errors in $L_{\infty}$ norm and the corresponding convergence rates.

\begin{table}[tbhp]
\tabcolsep=6pt
\caption{Errors and convergence rates with $\epsilon=10^{-9}$ in $L_{2}$ norm for  Problem 1.} \centering
\begin{tabular}{|c|c|c|cc|cc|cc|c|}
    \hline
    \multirow{2}{*}{Mesh}& \multirow{2}{*}{Iters} & \multirow{2}{*}{RRe} & \multicolumn{2}{c|}{$||U_h-u||_{2}$} & \multicolumn{2}{c|}{$||\widetilde{U}_h-u||_{2}$} & \multicolumn{2}{c|}{$||W_h-U_h||_{2}$}  &  \multirow{2}{*}{$r_h$} \\
    \cline{4-9}&   &   & Error &   Order &   Error &   Order  &   Error &   Order & \\
\hline
$  32\times  32\times  32$ &  8 & $2.24(-10)$ &$1.42(-4)$ &      &  $1.96(-7) $ &       &  $2.54(-5)$ &      &   $1.79(-1)$   \\
$  64\times  64\times  64$ &  9 & $2.68(-10)$ &$3.55(-5)$ & 2.00 &  $1.24(-8) $ & 3.98  &  $3.18(-6)$ & 2.99 &   $8.96(-2)$   \\
$ 128\times 128\times 128$ & 16 & $7.56(-10)$ &$8.87(-6)$ & 2.00 &  $7.83(-10)$ & 3.99  &  $3.99(-7)$ & 3.00 &   $4.50(-2)$   \\
$ 256\times 256\times 256$ & 78 & $8.98(-10)$ &$2.22(-6)$ & 2.00 &  $4.98(-11)$ & 3.97  &  $4.99(-8)$ & 3.00 &   $2.25(-2)$   \\
$ 512\times 512\times 512$ &  3 & $7.13(-10)$ &$5.55(-7)$ & 2.00 &  $3.06(-11)$ & 0.70  &  $6.25(-9)$ & 3.00 &   $1.13(-2)$   \\
\hline
\end{tabular}\label{table1}
\end{table}

\begin{table}[tbhp]
\tabcolsep=6pt
\caption{Errors and convergence rates with $\epsilon=10^{-9}$ in $L_{\infty}$ norm for Problem 1.} \centering
\begin{tabular}{|c|c|c|cc|cc|cc|}
    \hline
    \multirow{2}{*}{Mesh}& \multirow{2}{*}{Iters} & \multirow{2}{*}{RRe} & \multicolumn{2}{c|}{$||U_h-u||_{\infty}$} & \multicolumn{2}{c|}{$||\widetilde{U}_h-u||_{\infty}$} & \multicolumn{2}{c|}{$||W_h-U_h||_{\infty}$}  \\
    \cline{4-9}&   &   & Error &   Order &   Error &   Order  &   Error &   Order  \\
\hline
$  32\times  32\times  32$ &  8 & $2.24(-10)$ &$4.02(-4)$ &      &  $1.11(-6) $ &       &  $6.95(-5)$ &         \\
$  64\times  64\times  64$ &  9 & $2.68(-10)$ &$1.00(-4)$ & 2.00 &  $6.95(-8) $ & 4.00  &  $8.62(-6)$ & 3.01    \\
$ 128\times 128\times 128$ & 16 & $7.56(-10)$ &$2.51(-5)$ & 2.00 &  $4.35(-9) $ & 4.00  &  $1.07(-6)$ & 3.01  \\
$ 256\times 256\times 256$ & 78 & $8.98(-10)$ &$6.27(-6)$ & 2.00 &  $2.88(-10)$ & 3.92  &  $1.34(-7)$ & 3.00   \\
$ 512\times 512\times 512$ &  3 & $7.13(-10)$ &$1.57(-6)$ & 2.00 &  $1.15(-10)$ & 1.32  &  $1.67(-8)$ & 3.00    \\
\hline
\end{tabular}\label{table2}
\end{table}

As we can see that initial  guess $W_h$ is a third order approximation of the FE solution $U_h$, which validates our theoretical analysis in Section \ref{error}.
And the numerical solution $U_h$ reaches the full second accuracy for  both tolerances, while the extrapolated solution $\tilde{U}_h$ reaches fourth order accuracy on the coarse grids and starts to lose accuracy on fine grids.
This is due to the fact that the extrapolated solution $\tilde{U}_h$ is obtained from two second order numerical solutions $u_h$ and $u_{2h}$,
these two solutions must be extremely accurate in order to obtain a fourth order accurate solution $\tilde{U}_h$.
As the grid becomes finer, the tolerance needs to be more cruel.
 Thus, the extrapolated solution $\tilde{U}_h$ starts to lose convergence order when the grid is fine enough since a uniform tolerance is used in our ECMG algorithms.
 In addition, when the uniform tolerance becomes smaller, of course, the extrapolated solution $\tilde{U}_h$ reaches the full fourth order accuracy on more grids.  As shown in Table~\ref{table2}, on the grid $256\times256\times256$ and $512\times512\times512$ with the tolerance $\epsilon=10^{-9}$, the maximum error between the extrapolated solution $\tilde{U}_h$ and the exact solution $u$ reaches $O(10^{-10})$. The accuracy of $\tilde{U}_h$ is quite good, although the convergence order does not reach the full fourth order on the finest grid. If we want to obtain the fourth order accuracy on the finest grid, an even smaller  tolerance should be used.

We further point out that the number of iterations is reduced most significantly on the finest grid, which is particularly important when solving the large system. Below we provide a short illustration.
 Since we are using the stopping criteria which is related to the relative residual error in $L_2$ norm, in the following, we will show that the ratio $r_h$ measures how good $W_h$ approximates $U_h$ and  qualitatively reflects the number of iterations needed in the JCG solver with the initial guess $W_h$ in order to get the full second order accurate solution.
Since $||W_h-U_h||_{2}$ is third order convergent and $||U_h-u||_{2}$ is second order convergent, their ratio $r_h$ converges linearly to zero as the mesh is refined
(see the last column of Table~\ref{table1}).
Suppose the numerical solution on the previous grid (with mesh size $2h$) reaches the full second order convergence, then we have
\begin{align}
\parallel W_h-u\parallel_2&\leq \parallel W_h-U_h\parallel_2+\parallel U_h-u\parallel_2,\nonumber\\
&\leq (1+r_h)\parallel U_h-u\parallel_2,\nonumber\\
&\leq \frac{(1+r_h)}{4}\parallel U_{2h}-u\parallel_2.
\end{align}
If we use the initial guess $W_h$ directly as the numerical solution on the grid level $h$ without any iterations, the convergence order of $W_h$ is given by
\begin{align}
\textrm{order} =\log_2\frac{\parallel U_{2h}-u\parallel_2}{\parallel W_{h}-u\parallel_2}\geq 2-\log_2(1+r_h),
\end{align}
which means that the smaller the $r_h$ is, the more the convergence order $2-\log_2(1+r_h)$ approaches 2,
the less the iterations required to achieve the full second order accuracy are.
 Since $r_h$ decreases half when the grid is refined once, the convergence order $2-\log_2(1+r_h)$ will reach maximum on the finest grid.  Therefore, the number of iterations will be reduced most significantly on the finest grid. And this is particularly important when solving large linear systems.  For this example, we have $r_h\approx 0.00113$ on the finest grid $512\times 512\times 512$ for both tolerances.
 If using $W_h$ directly as the numerical solution on the finest grid, the convergent rate  is already greater than (or equal to)  $2-\log_2(1+r_h)=2-\log_2(1+0.00113)\approx 1.9838$, which is almost the full second order convergent. Therefore, only 3 iterations are needed to obtain the full second order convergent results, see the last row and second column of Table~\ref{table11}-\ref{table2}.

\begin{figure}[!tbp]
\centering
\subfigure[$k=0$]{
\label{Fig.sub.1}
\includegraphics[width=0.49\textwidth]{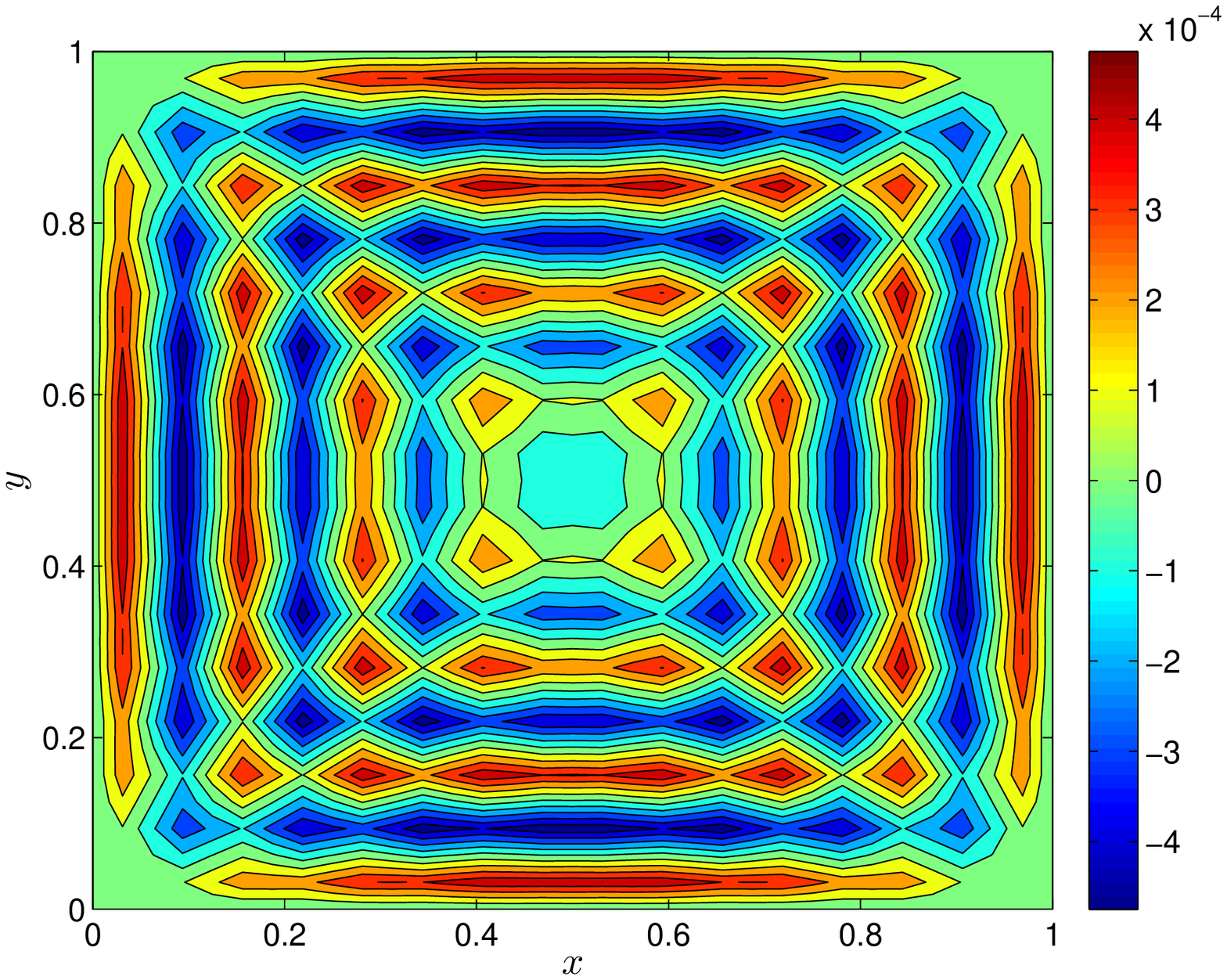}}
\subfigure[$k=1$]{
\label{Fig.sub.2}
\includegraphics[width=0.49\textwidth]{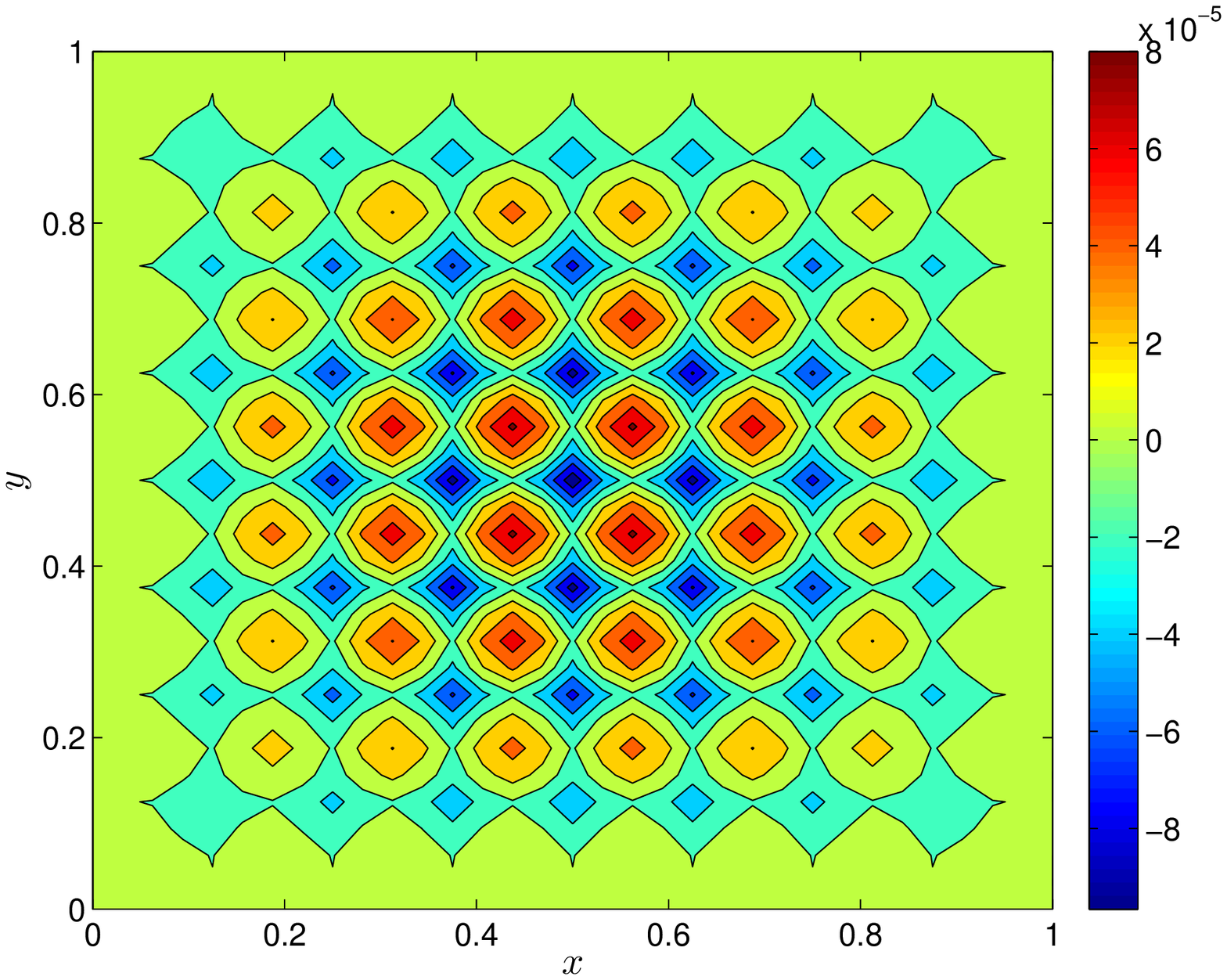}} \\
\subfigure[$k=3$]{
\label{Fig.sub.3}
\includegraphics[width=0.48\textwidth]{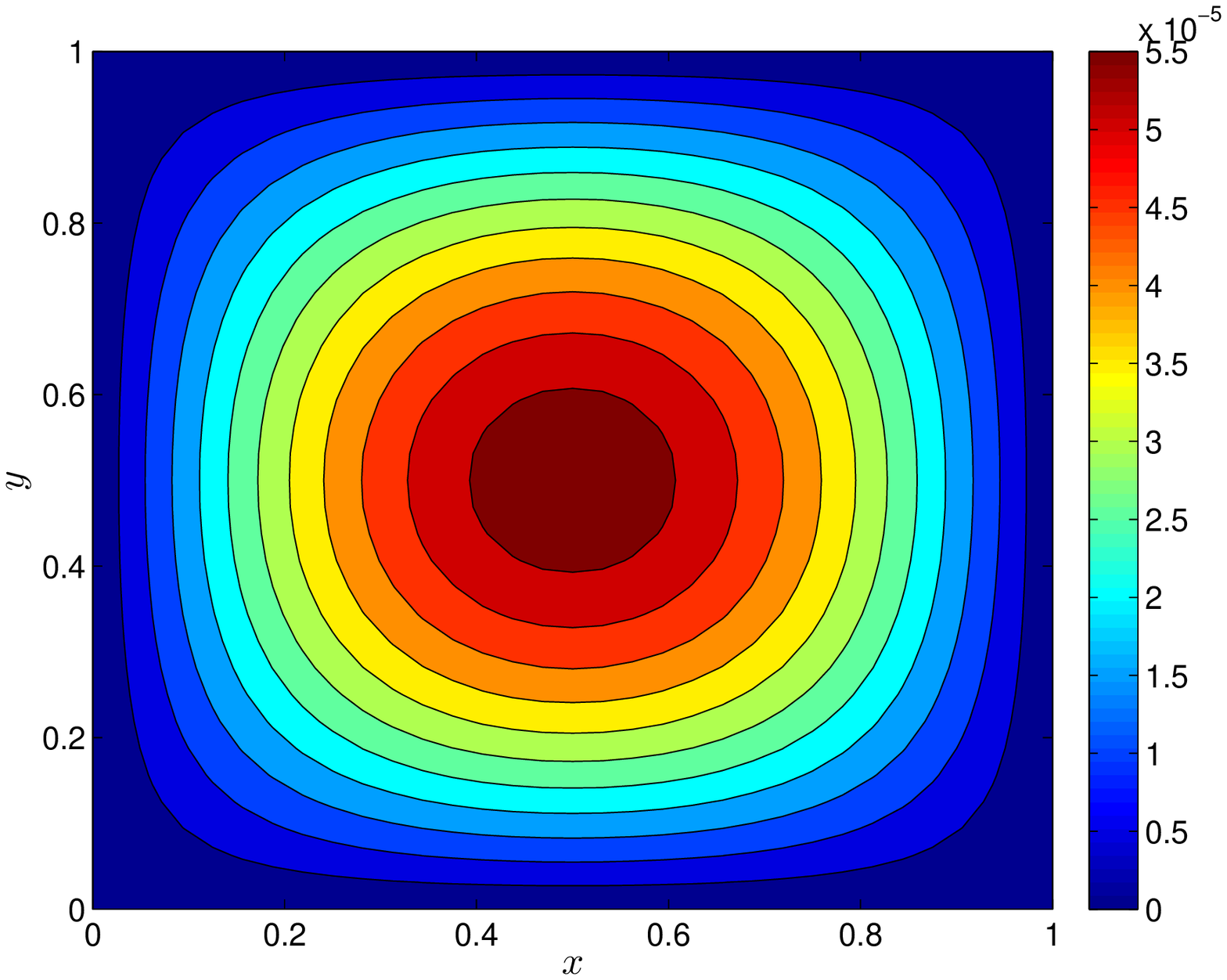}}
\subfigure[$k=8$]{
\label{Fig.sub.4}
\includegraphics[width=0.50\textwidth]{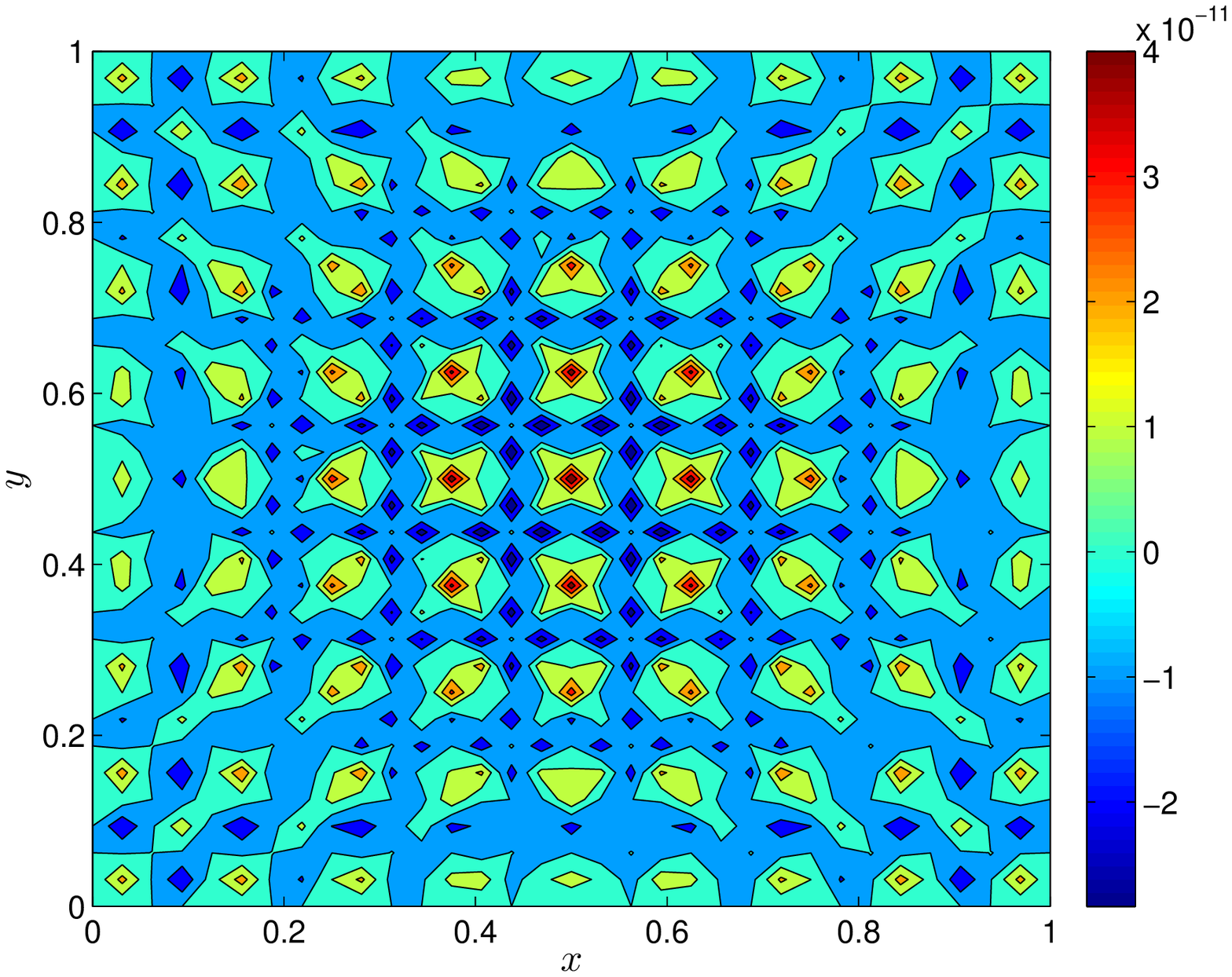}}\\
\caption{${U}_{1/32}^k-U_{1/32}$ on the plane $z=0.5$.}
\label{Fig.lable}
\end{figure}

Moreover, Figure~\ref{Fig.lable} presents the contour of the error (on the plane $z=0.5$) between exact FE solution $U_{1/32}$ and
 iterative solutions ${U}_{1/32}^k, k=0, 1, 3, 8$ of the JCG solver with $\epsilon=10^{-9}$, where ${U}_{1/32}^0=W_{1/32}$ is the initial guess. Since the error between the $W_{1/32}$ and $U_{1/32}$ is symmetric and oscillate,  the high frequency error components can be smoothed out easily. As we can see that the high-frequency oscillation is smoothed  out after only 3 iterations and the error has decreased one order of magnitude. In fact, the error between the  8-th iterative solution ${U}_{1/32}^8$ and $U_{1/32}$ is in the order of $O(10^{-11})$ while the relative residual $\frac{\parallel A_{1/32}{U}_{1/32}^8-f_{1/32}\parallel_2}{\parallel f_{1/32}\parallel_2}$ is $2.24\times 10^{-10}$,
which is less than $10^{-9}$ as shown in the first row and third column of Table~\ref{table1}. Thus, ${U}_{1/32}^8$ is actually the numerical solution $U_{1/32}$.

\begin{table}[!tbp]
\tabcolsep=6pt
\caption{Errors and convergence rates with $\epsilon=10^{-12}$ in $L_{2}$ norm for Problem 2.} \centering
\begin{tabular}{|c|c|c|cc|cc|cc|c|}
    \hline
    \multirow{2}{*}{Mesh}& \multirow{2}{*}{Iters} & \multirow{2}{*}{RRe} & \multicolumn{2}{c|}{$||U_h-u||_{2}$} & \multicolumn{2}{c|}{$||\widetilde{U}_h-u||_{2}$} & \multicolumn{2}{c|}{$||W_h-U_h||_{2}$}  &  \multirow{2}{*}{$r_h$} \\
    \cline{4-9}&   &   & Error &   Order &   Error &   Order  &   Error &   Order & \\
\hline
$ 40\times  16\times  20$ &  55 & $6.84(-13)$ &$2.97(-4)$ &      &  $4.81(-6)$ &       &  $5.93(-4)$ &     & $2.00(+0)$\\
$ 80\times  32\times  40$ &  81 & $8.84(-13)$ &$7.50(-5)$ & 1.99 &  $3.07(-7)$ &  3.97 &  $7.44(-5)$ & 2.99& $9.92(-1)$\\
$160\times  64\times  80$ & 137 & $9.47(-13)$ &$1.89(-5)$ & 1.99 &  $1.93(-8)$ &  3.99 &  $9.33(-6)$ & 3.00& $4.94(-1)$\\
$320\times 128\times 160$ & 136 & $9.82(-13)$ &$4.73(-5)$ & 2.00 &  $1.97(-9)$ &  3.30 &  $1.17(-6)$ & 3.00& $2.47(-1)$\\
$640\times 256\times 320$ &  12 & $9.95(-13)$ &$1.18(-6)$ & 2.00 &  $2.21(-9)$ & -0.17 &  $1.46(-7)$ & 3.00& $1.24(-1)$\\
\hline
\end{tabular}\label{table3}
\end{table}

\begin{table}[!tbp]
\tabcolsep=6pt
\caption{Errors and convergence rates with $\epsilon=10^{-12}$ in $L_{\infty}$ norm for Problem 2.} \centering
\begin{tabular}{|c|c|c|cc|cc|cc|}
    \hline
    \multirow{2}{*}{Mesh}& \multirow{2}{*}{Iters} & \multirow{2}{*}{RRe} & \multicolumn{2}{c|}{$||U_h-u||_{\infty}$} & \multicolumn{2}{c|}{$||\widetilde{U}_h-u||_{\infty}$} & \multicolumn{2}{c|}{$||W_h-U_h||_{\infty}$}\\
    \cline{4-9}&   &   & Error &   Order &   Error &   Order  &   Error &   Order  \\
\hline
$ 40\times  16\times  20$ &  55 & $6.84(-13)$ &$8.06(-4)$ &      &  $3.50(-5)$ &      &  $2.23(-3)$ &     \\
$ 80\times  32\times  40$ &  81 & $8.84(-13)$ &$2.02(-4)$ & 2.00 &  $2.50(-6)$ & 3.81 &  $2.78(-4)$ & 3.01\\
$160\times  64\times  80$ & 137 & $9.47(-13)$ &$5.04(-5)$ & 2.00 &  $1.68(-7)$ & 3.90 &  $3.47(-5)$ & 3.00\\
$320\times 128\times 160$ & 136 & $9.82(-13)$ &$1.26(-5)$ & 2.00 &  $1.06(-8)$ & 3.98 &  $4.34(-6)$ & 3.00\\
$640\times 256\times 320$ &  12 & $9.95(-13)$ &$3.14(-6)$ & 2.00 &  $6.68(-9)$ & 0.67 &  $5.43(-7)$ & 3.00\\
\hline
\end{tabular}\label{table4}
\end{table}

{\bf Problem 2.} The test Problem 2 can be written as
\begin{equation}\label{prob2}
\frac{\partial^2 u}{\partial x^2}+\frac{\partial^2 u}{\partial y^2}+\frac{\partial^2 u}{\partial z^2} =(1-2.5\pi^2)\exp(z)\sin(\frac{3\pi}{2} x)\sin(\frac{\pi}{2} y),\quad  \textrm{in } \Omega=[0,1]^3,
\end{equation}
where the boundary conditions are
\begin{equation}
  \frac{\partial u}{\partial n}(1,y,z)=\frac{\partial u}{\partial n}(x,1,z)=0
\end{equation}
and
\begin{equation}
  u(0,y,z)=u(x,0,z)=0,\quad u(x,y,0)=\sin(\frac{3\pi}{2} x)\sin(\frac{\pi}{2} y),\quad u(x,y,1)=e\sin(\frac{3\pi}{2} x)\sin(\frac{\pi}{2} y).
\end{equation}
The analytic solution of Eq.(\ref{prob3}) is
$$u(x,y,z)=\exp(z)\sin(\frac{3\pi}{2} x)\sin(\frac{\pi}{2} y),$$
where the exact solution is a sufficiently smooth function which changes more rapidly in the $x$ direction than in the $y$ and $z$ directions.

Since the solution has the fastest change in the $x$ direction and the slowest change in the $y$ direction, we use the coarsest grid $10\times4\times5$
in the ECMG algorithm.
Table~\ref{table3} and Table~\ref{table4} list the numerical data  obtained by ECMG$_{jcg}$ using a  tolerance $\epsilon=10^{-12}$.
Again, initial  guess $W_h$
is a third order approximation of the FE solution $U_h$. And the numerical solution $U_h$ reaches the full second accuracy, while the extrapolated solution $\tilde{U}_h$ reaches fourth order accuracy but starts to loss accuracy on the last two fine grids since we are using a uniform tolerance $\epsilon=10^{-12}$ on each level of grid. Additionally, the ratio $r_h$ converges to zero with order one.  If using $W_h$ as the numerical solution on the finest grid $640\times 256\times 320$, the convergent rate  is already greater than (or equal to)  $2-\log_2(1+0.124)\approx 1.8314$. Thus, only 12 iterations are required to achieve the full second order accuracy, see the last row and second column of Table~\ref{table3} and table~\ref{table4}.

{\bf Problem 3.} Consider a singular solution $u \in H^{3-\varepsilon}(\Omega) $ ($\varepsilon$ is any positive constant) satisfying
\begin{equation}\label{prob3}
\left\{ \begin{aligned}
         -\Delta u &= \frac{33xyz}{4(x^2+y^2+z^2)^{7/4}},\quad  &\textrm{in } &\Omega=[0,1]^3,\\
          u&=g(x,y,z),  &\textrm{on } &\partial \Omega,
        \end{aligned} \right.
\end{equation}
where $g(x,y,z)$ is determined from the exact solution
$$u(x,y,z) = \frac{xyz}{(x^2+y^2+z^2)^{3/4}}.$$
The exact solution $u$ has a removing singularity at the origin and has only finite regularity in $H^{3-\varepsilon}$.

Once again, we  use 7 level of grids with the coarsest grid $8\times8\times8$ and a  tolerance $\epsilon=10^{-11}$. Table~\ref{table5} lists the numerical data in sense of $L_2$ norm starting from the third level of grid, i. e., $32\times32\times32$. The initial  guess $W_h$ is still a third order approximation of the FE solution $U_h$, see second last column in Table~\ref{table5}.  And the numerical solution $U_h$ reaches the full second accuracy. Since the exact solution only has finite regularity in $H^{3-\varepsilon}$,  the extrapolated solution
$\tilde{U}_h$ can only reach third order accuracy, rather than fourth order accuracy for smooth solutions.
 Moreover, $r_h$ converges linearly to zero. If using $W_h$ as the numerical solution on the finest grid $512\times 512\times 512$, the convergent rate is already greater than (or equal to)  $2-\log_2(1+0.0966)\approx 1.87$. Thus, only 9 iterations are needed to obtain the full second order accurate solution, see the last row and second column of Table~\ref{table5}.

Since extrapolation are  based on asymptotic error expansions of the FE solution,
from Table \ref{table5} it seems that our ECMG method is still effective for such singular problems ($u\in H^{3-\varepsilon}$), and extrapolation can also
help us to increase the order of convergence to 3. This is a surprising result, which would widen the scope of applicability of our method.

\begin{table}[!tbp]
\tabcolsep=5pt
\caption{Errors and convergence rates with $\epsilon=10^{-11}$ in $L_{2}$ norm for Problem 3.} \centering
\begin{tabular}{|c|c|c|cc|cc|cc|c|}
    \hline
    \multirow{2}{*}{Mesh}& \multirow{2}{*}{Iters} & \multirow{2}{*}{RRe} & \multicolumn{2}{c|}{$||U_h-u||_{2}$} & \multicolumn{2}{c|}{$||\widetilde{U}_h-u||_{2}$} & \multicolumn{2}{c|}{$||W_h-U_h||_{2}$}  &  \multirow{2}{*}{$r_h$} \\
    \cline{4-9}&   &   & Error &   Order &   Error &   Order  &   Error &   Order & \\
\hline
$  32\times  32\times  32$ & 53 & $8.33(-12)$ &$2.80(-5)$ &      &  $2.25(-6) $ &       &  $3.23(-5)$ &      & $ 1.15(+0)$\\
$  64\times  64\times  64$ & 74 & $9.95(-12)$ &$7.16(-6)$ & 1.97 &  $2.88(-7) $ & 2.97  &  $4.56(-6)$ & 2.83 & $ 6.37(-1)$\\
$ 128\times 128\times 128$ & 52 & $9.32(-12)$ &$1.81(-6)$ & 1.98 &  $3.65(-8) $ & 2.98  &  $6.25(-7)$ & 2.87 & $ 3.45(-1)$\\
$ 256\times 256\times 256$ & 22 & $8.52(-12)$ &$4.57(-7)$ & 1.99 &  $5.14(-9) $ & 2.83  &  $8.42(-8)$ & 2.89 & $ 1.84(-1)$\\
$ 512\times 512\times 512$ &  9 & $7.60(-12)$ &$1.16(-7)$ & 1.98 &  $2.41(-9) $ & 1.09  &  $1.12(-8)$ & 2.91 & $ 9.66(-2)$\\
\hline
\end{tabular}\label{table5}
\end{table}

\section{Conclusions}
In this paper, we developed a  new extrapolation cascadic multigrid method, i.e., $\textrm{ECMG}_{jcg}$,  for solving the 3D elliptic boundary value problems on rectangular domains.  The major advantage of our method is to use the Richardson extrapolation and tri-quadratic Serendipity interpolation techniques for two numerical solutions on two level of grids  to obtain a quite good initial guess for the iterative solution on the next finer grid, which greatly reduces the iteration numbers for JCG solver.  In addition, a relative residual tolerance introduced in this paper can be used to control the accuracy of the numerical solutions more conveniently, and by using two second order numerical solutions on two scale grids, the
fourth order extrapolated solution $\tilde{U}_h$ on the fine grid can be obtained cheaply and directly.
Moreover, numerical results show that $\textrm{ECMG}_{jcg}$ has much better efficiency compared to classical MG methods  and is particularly suitable for solving large scale problems.

Our method developed in this paper can be easily extended to solve other related equations, for examples, convection-diffusion equations or Helmholtz equations. Moreover, the FE discretization method can be replaced by some other high order methods, such as compact finite difference methods. We are currently investigating these extensions.

\section*{Acknowledgements}
The authors would like to thank Professor Chuanmiao Chen for his useful discussions. And the help of gammar from Ms. Dedria Davis
is also grateful.


\begin{thebibliography}{10}

\bibitem{Avdeev2005}
D.~B. Avdeev, Three-dimensional electromagnetic modelling and inversion from theory to application. Surv. Geophys., 26 (2005), pp. 767-799 .

\bibitem{Newman2014}
G.~A. Newman, A Review of high-performance computational strategies for modeling and imaging of electromagnetic induction data. Surv. Geophys., 35 (2014), pp. 85-100.

\bibitem{Koldan2014}
J. Koldan, V. Puzyrev, J. de la Puente, et al. Algebraic multigrid preconditioning within parallel finite-element solvers for 3-D electromagnetic modelling problems in geophysics. Geophys. J. Int., 197 (2014), pp. 1442-1458.

\bibitem{Briggs2000}
W.~L. Briggs, S.~F. McCormick, V.~E. Henson, A Multigrid Tutorial, second ed.,
  SIAM, Philadelphia, PA, 2000.

\bibitem{Trottenberg2001}
U.~Trottenberg, C.~W. Oosterlee, A.~Sch¨¹ller, Multigrid, Academic Press,
  London, 2001.

\bibitem{Schaffer1984}
S. Schaffer, High order multi-grid methods, Math. Comp., 43 (1984), pp. 89-115.


\bibitem{Gupta1997a}
M. M. Gupta, J. Kouatchou, J. Zhang, Comparison of second and fourth order discretizations for multigrid Poisson solvers, J. Comput. Phys., 132(2) 226-232 (1997).

\bibitem{Othman1999}
M. Othman, A. R. Abdullah, An efficient multigrid Poisson solver, Int. J. Comput. Math., 71 (1999), pp. 541-553.


\bibitem{Zhang1998}
J. Zhang, Fast and high accuracy multigrid solution of the three dimensional Poisson equation, J. Comput. Phys., 143 (1998), pp. 449-461.

\bibitem{Zhang2002}
J. Zhang, Multigrid method and fourth-order compact scheme for 2D Poisson equation with unequal mesh-size discretization, J. Comput. Phys., 179 (2002), pp. 170-179.

\bibitem{Wang2009}
Y. Wang, J. Zhang, Sixth order compact scheme combined with multigrid method
and extrapolation technique for 2D poisson equation, J. Comput. Phys., 228 (2009), pp. 137-146.

\bibitem{Ge2010}
Y. B. Ge, Multigrid method and fourth-order compact difference discretization
scheme with unequal meshsizes for 3D poisson equation, J. Comput. Phys., 229 (2010), pp. 6381-6391.

\bibitem{Erlangga2006}
Y. A. Erlangga, C. W. Oosterlee, C. Vuik,  A novel multigrid based preconditioner for heterogeneous Helmholtz problems, SIAM J. Sci. Comput., 27 (2006), pp. 1471-1492.

\bibitem{Elman2001}
H. C. Elman,  O. G. Ernst,  D. P. O'leary,  A multigrid method enhanced by Krylov subspace iteration for discrete Helmholtz equations, SIAM J. Sci. Comput.,  23 (2001), pp. 1291-1315.

\bibitem{Erlangga20062}
Y. A. Erlangga, C. Vuik,  C. W. Oosterlee,  Comparison of multigrid and incomplete LU shifted-Laplace preconditioners for the inhomogeneous Helmholtz equation, Appl. Num. Math., 56 (2006), pp. 648-666.

\bibitem{Kim2002}
S. Kim, S. Kim,  Multigrid simulation for high-frequency solutions of the Helmholtz problem in heterogeneous media,  SIAM J. Sci. Comput., 24 (2002), pp. 684-701.


\bibitem{Riyanti2007}
C. D. Riyanti,  A. Kononov, Y. A. Erlangga, A parallel multigrid-based preconditioner for the 3D heterogeneous high-frequency Helmholtz equation. J. Comput. Phys., 224 (2007), pp. 431-448.


\bibitem{Zhang20022}
J. Zhang, H. Sun, J. J. Zhao, High order compact scheme with multigrid local mesh refinement procedure for convection diffusion problems, Comput.
Method Appl. M., 191 (2002), pp. 4661-4674.


\bibitem{Ge2011}
Y. B.  Ge, F. J. Cao, Multigrid method based on the transformation-free HOC scheme on
nonuniform grids for 2D convection diffusion problems, J. Comput. Phys.,  230 (2011), pp. 4051-4070.
%

\bibitem{Zeeuw1995}
P. M. De Zeeuw, El J. Van Asselt, The convergence rate of multi-level algorithms applied to the convection-diffusion equation. SIAM J. Sci. Stat. Comput., 6 (1985), pp. 492-503.

\bibitem{Gupta2000}
M. M. Gupta, J. Zhang, High accuracy multigrid solution of the 3D convection¨Cdiffusion equation. Appl. Math. Comput., 113 (2000), pp. 249-274.

\bibitem{Gupta1997}
M. M. Gupta, J. Kouatchou, J. Zhang, A compact multigrid solver for convection-diffusion equations. J. Comput. Phys., 132 (1997), pp. 123-129.

\bibitem{Bornemann1996}
F.~A. Bornemann, P.~Deuflhard, The cascadic multigrid method for elliptic
  problems, Numer. Math., 75 (1996), pp. 135-152.

\bibitem{Shaidurov1996}
V.~Shaidurov, Some estimates of the rate of convergence for the cascadic
  conjugate-gradient method, Comput. Math. Appl., 31 (1996), pp. 161-171.

\bibitem{Braess1999}
D.~Braess, W.~Dahmen, A cascadic multigrid algorithm for the stokes equations,
  Numer. Math., 82 (1999), pp. 179-191.

\bibitem{Timmermann2000}
G.~Timmermann, A cascadic multigrid algorithm for semilinear elliptic problems,
  Numer. Math., 86 (2000), pp. 717-731.

\bibitem{Shaidurov2000a}
V.~Shaidurov, L.~Tobiska, The convergence of the cascadic conjugate-gradient
  method applied to elliptic problems in domains with re-entrant corners,
  Math. Comput., 69 (2000), pp. 501-20.

\bibitem{Shaidurov2000}
V.~V. Shaidurov, G.~Timmermann, A cascadic multigrid algorithm for semilinear
  indefinite elliptic problems, Computing, 64 (2000), pp. 349-366.

\bibitem{Shi2000}
Z.~C. Shi, X.~J. Xu, Cascadic multigrid for parabolic problems, J. Comput. Math., 18 (2000), pp. 551-560.

\bibitem{Shi2001}
Z.~C. Shi, X.~J. Xu, A new cascadic multigrid, Sci. China Ser. A-Math., 44 (2001), pp. 21-30.

\bibitem{Braess2002}
D.~Braess, P.~Deuflhard, K.~Lipnikov, A subspace cascadic multigrid method for
  mortar elements, Computing, 69 (2002), pp. 205-225.

\bibitem{Stevenson2002}
R.~Stevenson, Nonconforming finite elements and the cascadic multi-grid method,
  Numer. Math., 91 (2002), pp. 351-387.

\bibitem{Bi2004}
C.~J. Bi, L.~K. Li, Cascadic multigrid method for isoparametric finite element
  with numerical integration, J. Comput. Math., 22   (2004), pp. 123-136.

\bibitem{Shi2004}
Z.~C. Shi, X.~J. Xu, H.~Y. Man, Cascadic multigrid for finite volume methods
  for elliptic problems, J. Comput. Math., 22 (2004), pp.  905-920.

\bibitem{Zhou2004}
S.~Z. Zhou, H.~X. Hu, On the convergence of a cascadic multigrid method for
  semilinear elliptic problem, Appl. Math. Comput., 159  2004), pp. 407-417.

\bibitem{Man2006}
F.~Y. Man, Z.~C. Shi, $P_1$-nonconforming quadrilateral finite volume element
  method and its cascadic multigrid algorithm for elliptic problems, J. Comput. Math., 24 (2006), pp. 59-80.

\bibitem{Shi2007}
Z.~C. Shi, X.~J. Xu, Y.~Q. Huang, Economical cascadic multigrid method (ECMG),
  Sci. China Ser. A-Math., 50 (2007), pp. 1765-1780.

\bibitem{Du2008}
Q.~Du, P.~B. Ming, Cascadic multigrid methods for parabolic problems, Sci. China Ser. A-Math., 51(8) 1415-1439 (2008).

\bibitem{Wang2008}
C.~Wang, Z.~Huang, L.~Li, Cascadic multigrid method for p-1-nonconforming
  quadrilateral element, J. Comput. Math., 16 (2008), pp. 237-248.

\bibitem{Xu2009}
X.~J. Xu, W.~B. Chen, Standard and economical cascadic multigrid methods for
  the mortar finite element methods, Numer. Math.-Theory Me., 2 (2009), pp. 180-201.

\bibitem{Yu2011}
H.~X. Yu, J.~P. Zeng, A cascadic multigrid method for a kind of semilinear
  elliptic problem, Numer. Algorithms, 58 (2011), pp. 143-162.

\bibitem{Chen2008}
C.~M. Chen, H.~L. Hu, Z.~Q. Xie, C.~L. Li, Analysis of extrapolation cascadic
  multigrid method (excmg), Sci. China Ser. A-Math., 51  (2008), pp. 1349-1360.

\bibitem{Hu2009}
H.~L. Hu, C.~M. Chen, Z.~Q. Xie, Extrapolation cascadic multigrid method
  (EXCMG)¡ªa new algorithm for solving large scale elliptic problems, Math.
  Numer. Sin., 31 (2009), pp. 261-274.

\bibitem{Chen2011}
C.~M. Chen, Z.~C. Shi, H.~L. Hu, On extrapolation cascadic multigrid method,
  J. Comput. Math., 29 (2011), pp. 684-697.

\bibitem{Marchuk1983}
G.~I. Marchuk, V.~V. Shaidurov, Difference Methods and Their Extrapolations,
  Springer-Verlag, New York, 1983.

\bibitem{Lin1983}
Q.~Lin, T.~Lu, S.-m. Shen, Maximum norm estimate, extrapolation and optimal
  point of stresses for finite element methods on strongly regular
  triangulation, J. Comput. Math., 1 (1983), pp. 376-383.

\bibitem{Lin1985}
Q.~Lin, J.-C. XU, Linear finite elements with high accuracy, J. Comput. Math., 3 (1985), pp. 115-133.

\bibitem{Blum1986}
H.~Blum, Q.~Lin, R.~Rannacher, Asymptotic error expansion and richardson
  extranpolation for linear finite elements, Numer. Math., 49 (1986), pp. 11-37.

\bibitem{Chen1989}
C.~M. Chen, Q.~Lin, Extrapolation of finite element approximation in a
  rectangular domain, J. Comput. Math., 7 (1989), pp. 227-233.

\bibitem{Fairweather2006}
G.~Fairweather, Q.~Lin, Y.~P. Lin, J.~P. Wang, S.~H. Zhang, Asymptotic
  expansions and richardson extrapolation of approximate solutions for second
  order elliptic problems on rectangular domains by mixed finite element
  methods, SIAM J. Numer. Anal., 44 (2006), pp. 1122-1149.

\bibitem{Asadzadeh2010}
  M.~Asadzadeh, A.~H. Schatz, W.~Wendland,
  Asymptotic Error Expansions for the Finite Element Method for Second Order Elliptic Problems in $R^N, N\geq2$. I: Local Interior Expansions,
  SIAM J. Numer. Anal., 48  (2010), pp. 2000-2017.

\bibitem{Liu2013}
G.~R. Liu, S.~S. Quek, The finite element method: a practical course,
  Butterworth-Heinemann, 2013.

\bibitem{Chen2009}
C.~M. Chen, H.~L. Hu, Z.~Q, Xie, et al. $L_2$-error of extrapolation cascadic multigrid (EXCMG), Acta Math. Sci., 29 (2009), pp. 539-551.

\end{thebibliography}
\end{document}